\documentclass{ametsocV6.1}

\usepackage{amsfonts}
\usepackage{amssymb}
\usepackage{amsmath}
\usepackage{graphicx}
\usepackage{setspace}
\usepackage{natbib}
\usepackage{xcolor}
\usepackage{subfig}
\definecolor{newcolor}{rgb}{.8,.349,.1}

\newcommand{\gray}[1]{\textcolor{gray}{#1}}

\newcommand{\mat}[1]{\boldsymbol{\mathsf{#1}}}
\renewcommand{\vec}[1]{\boldsymbol{\mathrm{#1}}}
\newcommand{\R}[1]{\mathbb{R}^{#1}}

\nolinenumbers

\title{Theoretical Advances in Current Estimation and Navigation from a Glider-Based Acoustic Doppler Current Profiler (ADCP)}

\authors{
    \nolinenumbers
    Jacob Stevens-Haas\aff{a} 
    \correspondingauthor{
        \nolinenumbers
        Jacob Stevens-Haas,  jmsh@uw.edu}, 
    Sarah E. Webster\aff{b}, 
    Aleksandr Aravkin\aff{a}
}

\affiliation{
    \nolinenumbers
    \aff{a}{Department of Applied Mathematics, University of Washington}.
    \aff{b}{Applied Physics Laboratory, University of Washington}
}

\abstract{
    \nolinenumbers
    We examine acoustic Doppler current profiler (ADCP) measurements from
    underwater gliders to determine glider position, glider velocity, and
    subsurface current.  ADCPs, however, do not directly observe the
    quantities of interest; instead, they measure the relative motion of
    the vehicle and the water column.  We examine the lineage of
    mathematical innovations that have previously been applied to this
    problem, discovering an unstated but incorrect assumption of independence.
    We reframe a recent method to form a joint probability model of current and
    vehicle navigation, which allows us to correct this assumption and extend
    the classic Kalman smoothing method.  Detailed simulations affirm the
    efficacy of our approach for computing estimates and their uncertainty.
    The joint model developed here sets the stage for future work to
    incorporate constraints, range measurements, and robust statistical
    modeling.
    \nolinenumbers
}

\begin{document}

\maketitle

%
%
%
\statement
The trade-offs between cost, power, and mission duration of underwater gliders allow them to measure underwater currents around the world.  We seek to improve the methods used by gliders to infer true current profiles and vehicle position from limited data.  Where previous methods applied smoothing techniques ad hoc, we start with a probability model for current and navigation from which we derive a more powerful and flexible smoothing method.  In simulation, our method performed well.  It was able to track within 75m a vehicle submerged and without GPS for three hours.  The flexible method allows experimenting with innovations from recent robust process model literature, including more difficult range measurements, and potentially real-time processing by low-powered gliders.
%
%
%

%

\vspace{.1in}
\section{Introduction}

    In order to map subsurface currents far away from fixed and mobile infrastructure, underwater gliders embark small, 1 MHz acoustic Doppler current profilers (ADCP).  These high-resolution profilers measure the relative velocity of a local slice of the water column and depend on both glider velocity and current. Additional measurements including GPS, acoustic long baseline beacons, and surface tenders as in \citet{Jakuba} can help with navigation.  Using this information along with process models of current variability and vehicle motion makes it possible to separate the relative measurement into its ground-referenced components: the subsurface ocean current profile and vessel navigation velocity.  The choice of prior process assumptions can have a significant impact on the quality of solutions, and the goal if this paper is to provide improved models for current and vehicle processes. 
    
    Previous efforts to develop a method for inferring currents and navigational data from ADCP measurements began with low-resolution sensors that imaged the entire water column.  The sensors typically descended from either a ship or fixed mooring.  Higher-resolution profilers had a limited range, and thus \citet{firing1990} developed the shear method to stitch together frames: the change in ADCP returns were averaged across overlapping depths to form the baroclinic part of the current profile.
    
    \citet{Visbeck2002} developed what is now the most widely used method, a linear least-squares model, commonly called the inversion method in reference to the pseudoinverse formed in solving such a model.  The inversion method originally involved two steps: one to measure depth-averaged current (DAC) and another to solve the matrix inverse problem.
    Methods derived from \citet{Visbeck2002} differ according to which terms are smoothed and how; whether the model includes other sensors such as bottom-tracking Doppler velocity log, inertial sensors, or terrain matching; whether the model treats currents as horizontally and temporally stable; and what measurement errors, such as sensor misalignment and ADCP bias, the method seeks to control.  The particulars of each method often reflect whether it aims to support underwater gliders or propeller-driven Autonomous Underwater Vehicles (AUVs).  \citet{Medagoda2016} contains a  literature review of both corpora, while we focus more specifically on methods for gliders.\footnote{
        The more complex methods generally aim to satisfy the needs of propeller-driven vehicles.  Such vehicles often conduct missions more concerned with precise navigation than current field modeling.  The methods also incorporate sensors that require more power and space than available on a typical underwater glider, such as the tactical-grade IMU of \citet{Medagoda2010}.
        
        Beyond the review in \citet{Medagoda2016}, valuable AUV-focused contributions have occurred in \citet{Arnold2018} and \citet{Medagoda2016a}, which test out different filters and a physics-based model, and \citet{Liu2021} which suggests a circular calibration dive prior to mission start and a high-pass filter in order to control IMU biases.
    }

    \subsection*{Contributions and Roadmap.}  
        We present an improved inversion method based on \citet{Jonker2019}, which frames the combined smoothing/ADCP inverse problem as a single joint optimization problem.\footnote{
            In the above taxonomy, their method, designed for underwater gliders, uses a process model smoother, minimal additional sensors, and a time-invariant, vertical current field.
        }  We start a step further back, with the stochastic model that gives rise to the optimization problem.  The theoretical starting point allows us to derive the proper smoothing covariance between the current and the vehicle's over-the-ground velocity and position.  We then derive the joint optimization problem from the posterior mode of vehicle and current states. The approach utilizes ADCP data, hydrodynamic model data, and GPS location fixes (at the start and, optionally, at the end of the dive). We demonstrate the success of the method in simulations without either a final GPS point, Doppler Velocity Log (DVL) seabed measurements, or Inertial Measurement Unit (IMU) measurements.
        
        The paper proceeds as follows.  We explain the inversion method and different ways authors have approached measurement and smoothing terms in Section \ref{lit_review}.  In Section \ref{methods} we describe the general state-space model, derive the covariance between current and vehicle process models, and describe the measurements we include in our model. We also derive higher order process model variants and discuss uncertainty quantification.  Section \ref{simulation} describes simulations and numerical experiments as well as their results. We conclude with a discussion and outline future work.

\section{Background}
    \label{lit_review}
    We begin this section with an explanation of the inversion method in \citet{Visbeck2002}, focusing on elements that this paper's antecedents have improved.  The basic method centers on the linear equation:
    \begin{align}
        \label{eqn:basic-linear}
        \vec b = \mat A \vec x,
    \end{align}
    where the vector $\vec b$ contains the set of ADCP measurements. The state vector, $\vec x$, represents the desired information as a stacked vector:
    \[
        \vec x := \left[\begin{matrix}
        \vec x_\upsilon\\
        \vec x_c
        \end{matrix}\right],
    \]
    with $\vec x_\upsilon$ as the vehicle velocity indexed by time and $\vec x_c$ as the  ocean currents indexed by depth.  The measurement matrix, $\mat A$,\footnote{The measurement matrix is also commonly referred to as the observation operator $\mat H$ in data assimilation problems.  However, Visbeck adds rows to the bottom of $\mat A$ and $\vec b$ to complete the data assimilation problem with smoothing terms.  Thus we use the more generic $\mat A\vec x=\vec b$ here.} selects the appropriate depth of current and time of vehicle motion for each ADCP measurement of relative velocity.  If $\vec x_\upsilon$ has entries at each time of $n$ measurements and $\vec x_c$ has entries for each depth of the $n$ measurements, then $\mat A\in \mathbb R^{2n}$ is column-rank deficient by $n$.  We would like the system to have full column rank in order to solve it as a least-squares problem using the pseudoinverse: 
    \[
        \vec x = \mat A^\dagger \vec b = (\mat A^T\mat A)^{-1}\mat A^T \vec b.
    \]
    We can remove $n-1$ degrees of freedom by smoothing adjacent currents and/or vehicle velocities.  We could instead remove $n-m$ degrees of freedom by binning nearby observations into the same $m < n$ depth/time bins.  If enough time bins have observations crossing depth bins or vice versa, we could remove all but one of the remaining degrees of freedom.  
    
    Nevertheless, even after binning and smoothing, we are left with at least one degree of freedom.  Every row of $\mat A$ in ADCP measurement or smoothing subtracts two values, which means that any $\vec x$ that satisfies equation \ref{eqn:basic-linear} would also be satisfied by $\vec x + \alpha \vec 1$ for some scalar $\alpha$.  Thus, regardless of binning or smoothing, the model requires absolute measurement of either current or vehicle kinematics, usually via GPS, in order for $\mat A$ to have sufficient column rank and the model be identifiable.
    
    \citet{Visbeck2002} discusses the effect of bin size choices for the depths considered in $\vec x_c$ as a trade off between  accuracy vs resolution.  However, as we will introduce a continuous model, we consider binning a choice between the numerical accuracy of the solution at the expense of a rigorous model of values between bin centers.  In addition, binning observations while smoothing helps the condition number of $\mat A^T\mat A$ by providing regular intervals between observations.\footnote{
        Locally, binning is the limit of smoothing as the smoothing penalty increases.  That is, smoothing allows a different current solution at each measurement depth, but forces those currents to be close together. Binning pools all nearby observations as measurements of the same bin depth.}

    \paragraph{Measurements.}
    Authors differ on how to treat the GPS measurement and depth-averaged current.  \citet{Visbeck2002} offers two mutually-exclusive measurement terms.  The first additional measurement term, that paper's expression (26), attempts to regularize the depth-average of current (DAC), or barotropic current, to be zero.  Thus, the regularizer implicitly tries to change $\vec x_c$ to model just {\it baroclinic} currents and shifts $\vec x_\upsilon$ by DAC.  We say ``tries'' because the choice of a regularizer rather than a constraint implies a belief that the average baroclinic current close zero, rather than exactly zero, as it is defined.  The paper then suggests that DAC or surface current could be calculated separately and added in to correct these values.  On the other hand, the alternative expression (30) in \citet{Visbeck2002} forces the integral of vehicle velocities to be equal to the displacement of the vehicle, as measured by GPS.  Such a term obviously requires $\vec x_\upsilon$ represent true, over-the-ground velocity and $\vec x_c$ represent the true current.

    \citet{Todd2017} identifies that the DAC regularization, expression (26) of~\citet{Visbeck2002}, can instead act as a measurement term if one calculates DAC before the inversion method. They calculate the DAC measurement as the difference between GPS measurement and dead reckoning: where the vehicle would have been if the hydrodynamic model's speed through the water is assumed accurate and current were absent.  He also includes a method to pre-calculate ADCP misalignment, pitch error biases, angle of attack, and sideslip.  
    Thus, while following the method of expression (26), Todd's two-step method still models vehicle over-the-ground velocity.  Although the reconciled DAC expression nominally would allow a GPS displacement measurement term like expression (30) of \citet{Visbeck2002}, the method already expends the GPS-measured displacement in calculating DAC.  It does, however, include a GPS-derived surface velocity measurement.

    Meanwhile, \citet{Jonker2019} reflects that one need not compose all the hydrodynamic model terms and GPS displacement into a single measurement of depth-averaged current; these terms may factor more properly  as individual measurements in a combined model.  Indeed, \citet{Visbeck2002} advocates the ease with which least-squares models can accommodate additional independent measurements.  Using GPS measurement directly precludes using it in preprocessing to calculate the ADCP and hydrodynamic model biases of \cite{Todd2011}.  However, measurement bias terms can be added into the state vector to try to achieve the same effect, as in \citet{Medagoda2015}.
    
    \paragraph{Smoothing.}
    Although \citet{Visbeck2002} bins the depths and times in state vector, the paper still offers a smoothing term for currents.  Smoothing terms in equation \eqref{eqn:basic-linear} are added to the bottom of $\mat A$ with zeros in the appropriate entry of $\vec b$.  The smoothing term in \cite{Visbeck2002} regularizes the second finite difference of current.  \citet{Todd2011} and \citet{Todd2017} follow the same method but smooth both the vehicle velocity and current.  Alternatively, \cite{Medagoda2010} introduces a Sparse Extended Information Filter in the context of AUVs. They follow it with different variants of Kalman filters in \cite{Medagoda2015}, \cite{Medagoda2016}, and \cite{Arnold2018}.  The Kalman filter derives from a Brownian motion model.  Such a model can give equivalent results in certain cases to finite difference smoothing, but finite difference smoothing does not derive as the discretization of a continuous stochastic process.  Their process model filter, in distinction to a smoother, calculates the optimal value from recent, real-time observations, as opposed to the greater power but greater expense of retaining the entire state history, as in the smoother of \citet{Jonker2019}.    

    The Kalman smoother requires modeling both the vehicle horizontal position as well as velocity, which simplifies the GPS measurement term.  The Kalman process variance term depends upon the interval between observations cubed, which means that binning at regular intervals helps limit the condition number of the method.  As a final note, \cite{Jonker2019} mentions that vehicle smoothing ought to properly apply to vehicle through-the-water velocity, not over-the-ground velocity.  However, because such smoothing would require modeling through-the-water velocity and consequent difficulties incorporating GPS measurement, \citet{Jonker2019} does not do so. 
    
    Overall, we contribute several important improvements to the method.  We explicitly model the problem as a maximum likelihood estimator, which, unlike standard inversion, is flexible enough to accommodate nonlinear measurements.  Our likelihood expression includes derivations of Kalman process terms for vehicle kinematics and  current.  The probability model lets us derive the correct covariances to smooth the vehicle's through-the-water velocity while still modeling the vehicle's over-the-ground velocity.  The process model approach also allows us to compare the smoothing between different process models.  We also eschew binning the data.
    
    
    

\section{Methods}
    \label{methods}
    Here, we specify the current and navigation random process terms and sensor measurement error terms.  The inversion method is equivalent to maximum likelihood estimation using ADCP measurement terms.  In the likelihood veiw, the addition of vehicle or current smoothing regularizers represents Bayesian prior knowledge of system behavior.  This representation is true whether the scientist explicitly proscribes a likelihood function or begins with the subsequent step, the matrix formulation of the solution method.
    
    Once we have specified independent likelihoods for each term, we combine them in a single optimization problem for the latent current and vehicle variables.  All of our terms are linear least squares, which means that the same solution method as prior authors (multiplication by a pseudoinverse) will work to solve our likelihood problem.
    
    Our contribution focuses on the process model priors.  It is largely compatible with other authors' innovations to handle measurement biases such as compass errors and ADCP misalignment.  We will ignore these innovations in the formulation of our model, as many of them deal with measurement terms or preprocessing, which will not be a factor in simulation and which we are not attempting to innovate.  However, we want to state the impact of ignoring these terms in the real world: We can expect an additional 17m of error for every degree of uncorrected compass error and every kilometer traveled.  Without an inertial measurement, we do not have a need for earth curvature terms; a local grid will not supply more than a 0.1\% error over the distances of a single dive between the 45th parallels.  As a final note, the assumption that preprocessing handles vehicle and geodetic orientation gives us all measurements in cardinal directions.  As such, the northward and eastward components of our model comprise two non-interacting problems.  When we specify a model term, we implicitly mean two terms: one for the northward component, and one for the eastward component.
    
    \subsection{Process Model}
    \label{xmodel}
    
    We use the vector $\vec x$ to describe the combined state space of the vehicle and current.  The state vector includes the vehicle's velocity and position on a local Cartesian grid and the water's velocity.  In avoiding binning observations, $\vec x$ has vehicle kinematic entries at all $n$ times an observation occurs and current entries at all $m$ depths an observation occurs.  We use a subscript to denote when selecting some subset of the state; e.g., $\vec x_\upsilon$ represents just the time-indexed vehicle components of the state (both position and velocity).  A complete list of this notation follows:
    \begin{enumerate}
        \item $\vec x_\upsilon \in \R{4n}$, vehicle components
        \item $\vec x_c \in \R{2m}$, current components
        \item $\vec x_{q}\in \R{2n}$, vehicle position components, one part of $\vec x_\upsilon$
        \item $\vec x_{\dot q} \in \R{2n}$, vehicle velocity components, the other part of $\vec x_\upsilon$
    \end{enumerate}

    As mentioned, the above vector space sizes include both northward and eastward components.

    While current is a deterministic process governed by fluid mechanics and vehicle motion is a deterministic process governed by basic kinematics and vehicle control, we need to specify how the processes behave absent much of the relevant knowledge.  In our initial model for $\mathrm{Pr}(\vec x)$, we assume the current varies by depth and vehicle velocity varies by time as samples of independent Brownian motions, which gives rise to standard Kalman smoothers in the likelihood expression.  That choice of randomness does not represent a physics or controls process that governs vehicle motion or current, but rather the structure of our uncertainty around these processes.  
    
    The Kalman smoother specifies the update matrix $\mat G_\upsilon$ and covariance matrix between increments of velocity and position, $\mat Q_\upsilon$.  It gives a negative log-likelihood term of:
    \begin{align}
        -\ell(\vec x_\upsilon) = \frac{1}{2\sigma_\upsilon^2}\|\mat G_\upsilon  \vec x_\upsilon\|_{\mat Q_\upsilon^{-1}}^2.
        \label{eqn:vehicle-process}
    \end{align}
    where $\sigma_\upsilon^2$ is a user-defined covariance scaling.
    
    Referring to time index $j$ for a coordinate with velocity $x_{\dot q}^j= W_{t_j}$ where $W_t$ is Brownian motion and position $x_q^j = \int_0^{t_j}W_r dr$, matrix $\mat G_\upsilon$ is responsible for creating mean-zero independent increments. As increments of Brownian motion are independent normal random variables, $\mat Q_\upsilon$  is block diagonal.	As a matrix equation, this expectation $\mathbb E$ takes the form
    \begin{align*}
        \left.\begin{matrix}
        	&\mathbb E\left[x^j_{\dot q} -x_{\dot q}^{j-1}|x_{\dot q}^{j-1}\right] &= 0\\
        	&\mathbb E\left[x^j_q-x_q^{j-1}- \Delta t \cdot x_{\dot q}^{j-1}|x_q^{j-1},x_{\dot q}^{j-1}\right]  &= 0
    	\end{matrix}\right. ~~\Rightarrow \mat G_\upsilon \vec x_\upsilon \sim \mathcal N(0, \mat Q_\upsilon).
	\end{align*}
	Diagonal blocks of $\mat G_\upsilon$ and $\mat Q_\upsilon$ have shape:
	\begin{align}
	    \label{eqn:basic-kalman-update}
    	\mat G_\upsilon^j = 
    	\left[\begin{matrix}
        	-1& 0 \\
        	-\Delta t_1 & -1 \\
    	\end{matrix}\right],~~~
    	\mat Q^j_\upsilon = 
    	\sigma_\upsilon^2\left[\begin{matrix}
        	\Delta t& \Delta t^2/2\\
        	\Delta t^2/2 & \Delta t^3/3
    	\end{matrix}\right],~~~
    \text{and}~~~~
    	\vec x_\upsilon = \left[\begin{matrix}
	    x_{\dot q}^1\\
	    x_q^1\\
	    \vdots\\
	    \end{matrix}\right],
	\end{align}
	and superdiagonal blocks of $\mat G_\upsilon$ are identity matrices.  Appendix A details the complete derivation of $\mat G_\upsilon$ and $\mat Q_\upsilon$.  
    
    The current process adds an analogous term, with three distinctions. For simplicity's sake, current varies by depth $s$ and not time $t$, because repeated measurements of the shear at the same depth are made over the course of less than a minute as the glider ascends or descends.  Secondly, in order to approximate the time/horizontal variance in the water column seen in different descent and ascent current profiles, we treat the ascent portion as merely continuing the dive.  Phrased differently, if we model a dive down to 750 m and back to the surface, we treat it as a single dive down to 1500 m.  Thus, currents at the same depth on ascent and descent are conditionally independent when conditioned on deeper currents.  Because of the straight paths and long dives that underwater gliders tend to fly, this simplification should not degrade results significantly.  Finally, the basic model only requires the water velocity and has no need for a current ``position.''  The negative log likelihood contribution for the current process mirrors the one for $\vec x_\upsilon$:
    \begin{align}
        -\ell(\vec x_{c}) = \frac{1}{2\sigma_c^2}\|\mat G_c \vec x_{c}\|_{\mat Q_c^{-1}}^2.
        \label{current_process}
    \end{align}
    A user-defined covariance scaling, $\sigma_c^2$, completes the term.

    We can compare these terms to the smoothers and regularizers employed in binning methods.  If intervals between observations remained constant, our method would regularize the first forwards difference of velocity.\footnote{However, the Kalman smoother also regularizes vehicle position according to its covariance with the unmodeled velocity between bins.} \citet{Visbeck2002} and \citet{Todd2011} instead regularize the 2\textsuperscript{nd} centered difference without integrating that error into position.  Their finite difference regularizer does not have a continuous stochastic representation akin to our model.

    \subsection{Measurement Terms}
    \label{zmodel}
    In addition to the priors responsible for smoothing the current and vehicle processes, we need measurement terms.  We use ADCP measurements of the water column's velocity at one depth relative to the vehicle over-the-ground velocity at the same time.  A hydrodynamic model estimates vehicle through-the-water velocity,  which is mathematically similar to an ADCP measurement as noted in \cite{Todd2017}, but allowed a different error variance.  Finally, a pair of GPS fixes at either the beginning of the dive or spanning it measure the change in vehicle position.
    
    We present the measurement error expressions below. $\mat B^q$ and $\mat B^{\dot q}$ select the vehicle  kinematics of the state vector at the appropriate measurement times, and $\mat B^c$ selects the currents of the state vector at the appropriate measurement depths.  Each term has a user-specified variance parameter, $\sigma^2$, representing the measurement error variance; one can set these parameters based on published sensor accuracy.  None of these error terms includes a modeled bias or covariance across measurements.  Again, we only write the explicit terms for either the eastward or northward components; in actuality, we use two of each term.  
    
     \begin{enumerate}
    	\item $\vec z_{adcp}$: ADCP data measured from the glider with error variance $\sigma^2_{adcp}$ (In our case, \cite{NortekADCP1000} describes the variance of the Nortek Signature 1000 ADCP as $1e-3 m^2/s^2$).  This measurement reflects the current at a certain depth relative to the true velocity of the glider: 
    	\begin{align}
    	\label{loss-adcp}
            -\ell(\vec z_{adcp}|\vec x)=\frac{1}{2\sigma_{adcp}^2}\|\vec z_{adcp}-(\mat B^{\dot q}_{adcp}\vec x_{\dot q}-\mat B^c_{adcp} \vec x_{c})\|^2.
        \end{align}
    	\item $\vec z_{ttw}$: Through-the-water estimated velocity of the glider with error variance $\sigma^2_{ttw}$. While this quantity nominally relies on calculations from pitch and depth rate measurements using a vehicle's hydrodynamic model, we treat it as direct measurement.  For real dives, models such as \cite{Eriksen2001} provide calculations based on pressure, buoyancy, and pitch.  That paper measured a vertical velocity error variance of $1.58 cm^2/s^2$.  Horizontal velocity variance is a function of angle of attack, but for a lift/drag ratio of 3:1, the modeled horizontal velocity variance would be around $14.2 cm^2/s^2$ 
    	\begin{align}
    	\label{loss-ttw}
            -\ell(\vec z_{ttw}|\vec x)=\frac{1}{2\sigma_{ttw}^2}\|\vec z_{ttw}-(\mat B^{\dot q}_{ttw} \vec x_{\dot q}-\mat B^{c}_{ttw}\vec x_c)\|^2 .
	    \end{align}
	    \item $\vec z_{gps}$: GPS position measurements at the surface before and after a dive with error variance $\sigma^2_{gps}$.  Practically, $\sigma^2_{gps}$ is usually in the $O(1)$ range.  As a representative value, \cite{Garmin} lists the GPS15xH to have 3-5 m accuracy 95\% of the time, which works out to a variance of $2.5-6.5 m^2$.
    	\begin{align}
    	    \label{eqn:loss-gps}
            -\ell(\vec z_{gps}|\vec x)=\frac{1}{2\sigma_{gps}^2}\|\vec z_{gps}-(\mat B^q_{gps} \vec x_q)\|^2.
	    \end{align}
    \end{enumerate}
    
    These terms compose the complete likelihood specification for the model:
    \begin{align}
        \label{eqn:basic-likelihood}
        -\ell(\vec x|\vec z) = -\ell(\vec x_{\upsilon})-\ell(\vec x_{c}) -\ell(\vec z_{p}|\vec x) -\ell(\vec z_{ttw}|\vec x) -\ell(\vec z_{gps}|\vec x).
    \end{align}
    
    While other authors model more robust measurement terms, we omit terms for compass error, ADCP bias due to misalignment, or hydrodynamic model bias due to sideslip and angle-of-attack errors aggravated by biofouling.
    
    \subsection{Extensions}
    \label{variations}
        Here we discuss the extensions and variations of our basic model. 
        
        \subsubsection{Higher-order Smoothing}
        \label{sec:higher-order}
        In Subsection~\ref{methods}.\ref{xmodel}, we described vehicle over-the-ground velocity and current as Brownian processes.  The vehicle's over the ground position then becomes the integral of Brownian motion.  Whereas previous authors have described the concomitant regularizer as a smoothing term, strictly speaking, the assumption of Brownian motion implies nonsmoothness.  So we can assume smoothness if we expand our state space model to include higher order terms: vehicle acceleration and the change in current with respect to depth.  These terms then become Brownian.  Their integrals, vehicle velocity and current velocity, become smooth.
        
        This decision changes the process update matrices, $\mat G_\upsilon$ and $\mat G_c$, and the covariance matrices, $\mat Q_\upsilon$ and $\mat Q_c$:
    	\begin{align}
        	\label{eqn:higher-order-kalman}
        	\mat G^j_\upsilon = 
        	\left[\begin{matrix}
            	-1& 0 & 0 & \\
            	-\Delta t_1 & -1 & 0 & 0 & \\
            	-\frac{\Delta {t_1}^2}{2}& -\Delta t_1 & -1 \\
        	\end{matrix}\right],
        	\qquad 
        	 \mat Q^j_\upsilon=\left[\begin{matrix}
            	{\Delta t} & \frac{\Delta t^2}{2} & \frac{\Delta t^3}{6}\\
            	\frac{\Delta t^2}{2} & \frac{\Delta t^3}{3} & \frac{\Delta t^4}{8}\\
            	\frac{\Delta t^3}{6} & \frac{\Delta t^4}{8} & \frac{\Delta t^5}{20}\\
        	\end{matrix}\right],
        	\qquad
        	\text{and}~~~~
        	\vec x_\upsilon = \left[\begin{matrix}
    	    x_{\ddot q}^1\\
    	    x_{\dot q}^1\\
    	    x_q^1\\
    	    \vdots\\
    	    \end{matrix}\right].
    	\end{align}
    	As before, $\mat G_\upsilon$ has superdiagonal identity matrix blocks.
        
        Now that our current process has more than one order, we use $\vec x_{\dot p}$ to refer specifically to the velocity components of the current process and $\vec x_{\ddot p}$ to refer to the change in velocity with respect to depth.  $\vec x_c$ then refers more generally to all the current process terms in the same way that $\vec x_\upsilon$ refers to components $\vec x_q$, $\vec x_{\dot q}$, and now $\vec x_{\ddot q}$, vehicle acceleration.  Appendix B contains more details.
    
        \subsubsection{Independence.}
        \label{sec:independence}
        Unfortunately, our likelihood expression \ref{eqn:basic-likelihood} implies an incorrect assumption.  Adding together $-\ell(\vec x_\upsilon)$ and $-\ell(\vec x_c) $ implies that the joint probability $\mathrm{Pr}(\vec x_\upsilon, \vec x_c) = \mathrm{Pr}(\vec x_\upsilon)\mathrm{Pr}(\vec x_c)$, which in turn implies independence of the vehicle's over-the-ground velocity from current.  As far as we can tell, all previous authors who modeled over-the-ground vehicle kinematics and additively smoothed both current and vehicle kinematics made the same assumption \textit{de facto}.  Instead, the glider's {\it through-the-water} velocity is independent of current;  its over-the-ground velocity obviously depends upon current, and this dependence becomes more significant in larger currents: 
        \begin{align}
        \label{relative}
            \vec x_{\dot q} = \vec x_{\dot p} + \vec x_{\dot r},
        \end{align}
        
        where $\vec x_{\dot r}$ is the relative velocity of the vehicle through the water.
        
        To correct our likelihood formula, we have two options.  (1) We can switch from modeling $\vec x_q$ and $\vec x_{\dot q}$ to modeling $\vec x_r$ and $\vec x_{\dot r}$. (2) Alternatively, we can use the conditional formula, $\mathrm{Pr}(\vec x_\upsilon, \vec x_c) = \mathrm{Pr}(\vec x_\upsilon|\vec x_c)\mathrm{Pr}(\vec x_c)$ in order to keep modeling $\vec x_\upsilon$.  Both involve similar computations.
        
        \paragraph{Modeling $\vec x_r$.}  In the former variant, we model the vehicle's kinematics relative to the water column.  However, if $\vec x_{\dot r}$ represents the vehicle's through-the-water velocity, then $\vec x_r$ represents the vehicle's position relative to the water.  In order to incorporate GPS measurements, we must calculate the true, over the ground position of the vehicle by integrating equation \eqref{relative}.  To effect that calculation, we initially let the current velocity field, $x_{\dot p}$ vary as a function of depth and time.  We integrate equation \eqref{relative}:
    	\begin{align}
    		x_q(T) - x_q(0)&= \int_{0}^{T} x_{\dot q}(t) dt = \int_{0}^{T} x_{\dot r}(t) + x_{\dot p}(s(t),t) dt.
    	\end{align}
    	We change variables using $\dot s=ds/dt$ in the second integral term:
    	\begin{align}
        	\label{change-of-variables}
    		&= x_r(T)-x_r(0) + \int_{s(0)}^{s(T)}x_{\dot p}(s, t(s))\frac{1}{\dot s} ds\\
    		&= x_r(T)-x_r(0) + x_p(s(T)) - x_p(s(0)).
    	\end{align}
    	Here, we recognize $x_p(s)$ as the position of a test particle drifting in the current but changing depths to experience the same currents as the vehicle.  The vehicle's over-the-ground position change would then equal the sum of the test particle's position and the vehicle's through-the-water position.  Our state vector $\vec x_c$ now comprises both $\vec x_{\dot p}$ and $\vec x_p$.
    	
    	Practically, the additional current process order manifests as a different $\mat G_c$ and $\mat Q_c$ similar to that of vehicle kinematics in the base method but with the $\dot s^{-1}$ term from equation \ref{change-of-variables}.  Appendix C contains more details.
    	
    	In addition to those smoothing terms, the $\mat B_{adcp}$ and $\mat B_{ttw}$ matrices adjust to select the proper velocities for their respective measurements.  Moreover, the GPS measurement term, expression \ref{eqn:loss-gps}, now requires both a $\mat B_{gps}^p$ and a $\mat B_{gps}^r$ to add up the true position change of the vehicle.  
    
        \paragraph{Covariance between $\mathrm{Pr}(\vec x_\upsilon)$ and $\mathrm{Pr}(\vec x_c)$.}  If instead we model vehicle kinematics over the ground (as in the basic method), we must derive our maximum likelihood estimate from:
        \begin{align}
        \label{eqn-conditional}
            \mathrm{Pr}(\vec x_q, \vec x_{\dot q}, \vec x_c) = \mathrm{Pr}(\vec x_q, \vec x_{\dot q} | \vec x_c) \mathrm{Pr}(\vec x_c)\\
            \ell(\vec x_q, \vec x_{\dot q}, \vec x_c) = \ell(\vec x_q, \vec x_{\dot q} | \vec x_c)+ \ell(\vec x_c).
        \end{align}
        We achieve the conditional likelihood expression by reusing equation \eqref{change-of-variables} to get the joint probability, then conditioning on $\vec x_c$.  This gives rise to a different $\mat G_\upsilon$ and $\mat Q_\upsilon$.  $\mat G_\upsilon$ is the same as in the base case, but expression \ref{eqn:vehicle-process} needs to accomodate the conditioning on $\vec x_c$.  We introduce $\widetilde{\mat G}_\upsilon$ to capture the effect of $\vec x_c$ on $\vec x_\upsilon$ in the likelihood expression
        
        \begin{align}
            \ell(\vec x_q, \vec x_{\dot q} |\vec x_c) = \frac{1}{2\sigma_\upsilon^2}\|\mat G_\upsilon  \vec x_\upsilon+\widetilde{\mat G}_\upsilon  \vec x_c\|_{\mat Q_\upsilon^{-1}}^2,
        \end{align}
        where $\widetilde{\mat G}_\upsilon$ has 2x1 blocks.  Letting $k(j)$ refer to the index of current experienced by the vehicle at time $t_j$ (that is, $s_{k(j)} = s(t_j)$) and $M$ refer to the total number of timesteps, we have:
        \begin{align}
        \label{eqn:conditional-gq}
            \widetilde{\mat G}_\upsilon^{j, M+k(j)} =\left[\begin{matrix}
                \Delta x_c \\
                \frac{\Delta s}{2\dot s}\Delta x_c
            \end{matrix}\right]
            \qquad
            \widetilde{\mat G}_\upsilon^{j, M+k(j+1)} =\left[\begin{matrix}
                -\Delta x_c \\
                -\frac{\Delta s}{2\dot s}\Delta x_c
            \end{matrix}\right]
            \qquad
            \mat Q^j_\upsilon = 
                \left[\begin{matrix}
                    \Delta t  & \frac{\Delta t^2}{2}\\
                    \frac{\Delta t^2}{2} & \frac{\Delta t^3}{3}+\frac{\sigma^2_c}{\sigma^2_\upsilon}\frac{\Delta s^3}{12\dot s^2}
                \end{matrix}\right].
        \end{align}
        
        Appendix D provides a step-by-step derivation, including the matrices when $\vec x_\upsilon$ and $\vec x_c$ smoothed an additional order or when conditioning on depths in between $s_{k(j)}$ and $s_{k(j+1)}$.

    \subsection{Optimization and Uncertainty Quantification}
    \label{uncertainty}
        
        The likelihood approach of expression \eqref{eqn:basic-likelihood} corresponds to a Bayesian interpretation of the loss. Our prior belief of the Brownian process terms, combined with the Gaussian measurement likelihoods, results in a combined posterior and its estimator that are Gaussian random variables.  Given expression \eqref{eqn:basic-likelihood}, one can represent the posterior estimator of $\vec x$ as:
        \begin{align}
            \label{eqn:complete-likelihood}
            \hat{\vec x} = \arg\min_{\vec x}&-\ell(\vec x) \\
            =\arg\min_{\vec x}
            &\frac{1}{2\sigma^2_\upsilon}\|\mat G_\upsilon \vec x_\upsilon\|^{2}_{\mat Q_\upsilon^{-1}} + \frac{1}{2\sigma^2_c}\|\mat G_c\vec x_c\|^2_{{\mat Q_c}^{-1}} + \nonumber \\
            &\frac{1}{2\sigma^2_{adcp}}\left\|\vec z_{adcp}- (\mat B^{\dot q}_{adcp}\vec x_{\dot q} - \mat B^{\dot p}_{adcp} \vec x_{\dot p})\right\|^{2} + \\
            &\frac{1}{2\sigma^2_{ttw}}\left\|\vec z_{ttw}- (\mat B^{\dot q}_{ttw}\vec x_{\dot q} - \mat B^{\dot p}_{ttw} \vec x_{\dot p})\right\|^{2}+ \nonumber \\
            &\frac{1}{2\sigma^2_{gps}}\left\|\vec z_{gps}- (\mat B^{q}_{gps}\vec x_{q} - \mat B^{p}_{gps} \vec x_{ p})\right\|^{2}. \nonumber 
        \end{align}
        One can compute the estimator in closed form, using the normal equations.  Knowing that any subscript on $\vec x$ practically represents matrix multiplication and that both $\mat Q^{-1}$, as positive definite matrices, factor to $LL^T$, we can rewrite every term in equation \ref{eqn:complete-likelihood} as 
        \begin{align}
            \frac{1}{2\sigma^2_i}\left\|\mat A_ix-\vec b_i\right\|^2.
        \end{align}
        And stacking all $\mat A_i/\sigma^2_i$ into matrix $\mat A$ and all $\vec b_i/\sigma^2_i$ into vector $\vec b$, we can write the optimal value as:
        \begin{align}
            \nabla \ell(\hat{\vec x}) &= 0\\
            &=\mat A^T \mat A \hat{\vec x} - \mat A^T \vec b\\
            \hat{\vec x} &= \mat A^\dagger \vec b = (\mat A^T\mat A)^{-1}\mat A^T\vec b.
        \end{align}
        
        Moreover, the variance of the estimator $\vec {\hat x}$ has an explicit solution:
        \begin{align}
            \textrm{var}(\hat{\vec x}) = \mathbb E [(\hat{\vec x} - \vec x )(\hat{\vec x} - \vec x )^T] = (\mat A^T\mat A)^{-1}.
        \end{align}
        
    \subsection{Simple Comparisons: Dead Reckoning, Corrected for Depth-Averaged Current}
        When we run trials and develop solutions via the inversion method, we want to be able to compare our results to those generated by a more naive method.  Dead Reckoning (DR) is the basic navigation method where, in the absence of external positioning information and current information, a vessel advances its position purely based upon inaccurate velocity information.  In this case, we dead reckon the vehicle's position based solely upon the relative velocity calculated by the hydrodynamic model.  After resurfacing, a new GPS fix can tell us how far off our dead reckoned position is, and assuming all error is due to currents, we can calculate the depth-averaged current (DAC).  We can then adjust our vehicle velocities at each point by DAC to generate an overall navigational track. Glider's have routinely used the difference between the end-of-dive GPS and dead reckoning as an estimate of the true depth-averaged current, accurate to $\sim1$ cm/s in order to successfully aid in navigation. 

\section{Simulation Results}
\label{simulation}
    \subsection{Simulation Setup}
        We conducted simulations to compare our method and its variants.  Each trial includes a variant with GPS points at both ends of the dive and an identical variant with two GPS points at the beginning and none at the end.  This second variant behaves similarly to how \citet{Todd2017} includes initial surface drift velocity, as measured by GPS, as a term in their least-squares solution.  However, we leave the observations as positional in nature, which allows the Kalman machinery to manage the covariance between positions and velocities.  As discussed later in the conclusion, this ability to swap in different measurements shows the flexibility of the probability model framework.
        
        The dive simulation parameters reflect a randomized current profile and vehicle navigation path, with max current and max propulsion speed similar to those seen in Seaglider experiments in the arctic in \citet{Jonker2019}.  We made compromises for ease of simulation, such as separate, uncorrelated sinusoids for current and vehicle propulsion on ascent and descent; this manifests as an apparent vehicle turn at apogee.  Moreover, sinusoids are a common simulation for Kalman Smoothing demonstrations.  Because the maximum current and propulsion speeds are normally distributed, they occasionally exceed those likely in reality.  We believe that model success on these profiles will translate well to more modest profiles.  Table \ref{tab:sim-params} lists simulation parameters.
        \begin{table}[h]
            \centering
            \begin{tabular}{c|c}
                Parameter &  Value\\
                \hline
                Duration & 3 hrs \\
                Depth & 750 m \\
                Hydrodynamic model measurements & 500 \\
                ADCP measurements & 450 \\
                ADCP bins & 4 \\
                ADCP range & 12m, upwards-facing \\
                Current & piecewise sinusoid with amplitude $\sim N(0, .3)$ kts \\
                Vehicle TTW velocity & piecewise sinusoid with amplitude $\sim N(0, .4)$ kts \\
                Hydrodynamic model noise & 1e-2 m/s \\
                ADCP measurement noise & 1e-2 m/s\\
                GPS measurement noise & 1 m
            \end{tabular}
            \caption{Simulation Parameters, reflecting likely dive of a Seaglider with Nortek ADCP}
            \label{tab:sim-params}
        \end{table}

    \subsection{Parameters and Grid-Search}
        While published information about physical devices can inform setting the measurement variance parameters $\sigma^2$, the process variance parameters for smoothing represent prior knowledge of vehicle and current dynamics.  We wish to understand how precisely our model's success depends upon knowing the exact process model variances, parameters we can only ever assume imprecisely\footnote{Formally, the problem of simultaneous state and variance estimation for Kalman smoothing diverges.  Some attempts have been made for Kalman filter variance identification using the Normalized Innovation, Squared metric in \cite{Cai2019}.  Unfortunately, cases exist where the likelihood function of $\sigma^2_\upsilon$ and $\sigma^2_c$ are indeterminate.  Thus, rather than finding a particular best variance parameter, we are interested in the size of $\sigma_\upsilon^2$, $\sigma_c^2$ parameter space that gives a reasonable solution}.  Thus, in each simulation, we run a parameter search on the variance coefficients for the process terms, $\sigma_\upsilon^2$ and $\sigma_c^2$, and evaluate the effect on accuracy metrics.  The measurement error terms in modeling, $\sigma_{adcp}^2$, $\sigma_{ttw}^2$, and $\sigma_{gps}^2$, are given the same value used in simulation.

    \subsection{Results}
\begin{figure}
            \centering
            \includegraphics[width=32pc]{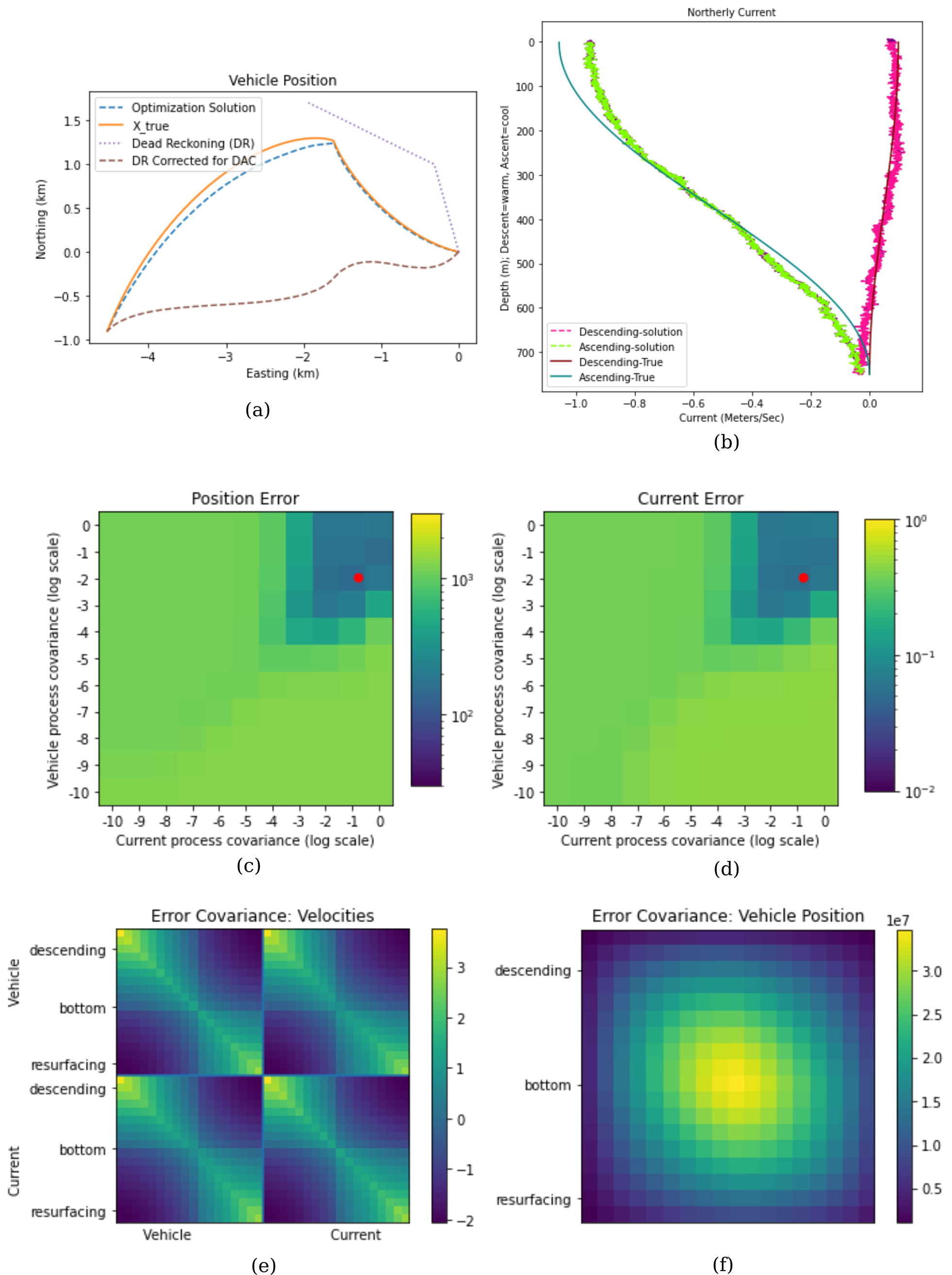}
            \caption{Visualizations of evaluation criteria for solution method 1 with both starting and ending GPS measurements.  (a) Navigation solution, (b) Current solution (only northward shown; eastward error shows similar visual fit), (c) Navigation error across process covariance parameter search, with minimizer marked by red dot (d) Current error across process covariance parameter search, with minimizer marked by red dot(e) Vehicle and current velocity estimator uncertainty, and (f) Vehicle position estimator uncertainty.  The parameter search plots average twenty simulations.  The rest of the plots represent a single simulation from among those twenty.  See appendix F for all plots from all trials, organized for comparison across trials.}
            \label{fig:all}
        \end{figure}

        We compare four methods based on different process models: (1) the basic one described in section \ref{methods}.\ref{xmodel}, (2) the higher-order one described in section \ref{methods}.\ref{variations}.\ref{sec:higher-order}, (3) the process with proper vehicle-current covariance described in section \ref{methods}.\ref{variations}.\ref{sec:independence}, and (4)  a combination of methods 3 and 4: a proper covariance method that smooths to higher order.  While section \ref{methods}.\ref{variations}.\ref{sec:independence} argues that modeling through-the-water kinematics is mathematically equivalent to modeling true kinematics with appropriate vehicle-current covariance, in practice, we focused on developing the solution that was easier to code and verify.  Thus, we only display results for the latter method.  In comparing methods, we use as criteria: (a) the horizontal position accuracy during the dive at optimal process parameters\footnote{
            When the process parameters that give rise to the most accurate navigation estimates and those that give rise to the most accurate current estimates disagree, we choose the parameters that lead to the most accurate current estimation as optimal.
        }, (b) the current accuracy during the dive at optimal process parameters, (c) the breadth of process parameters that achieve satisfactory results, and (d) performance on simulations including only two closely-placed GPS fixes before the dive and no GPS fix at the end.  We evaluate root-mean-square-error (RMSE) across the entire vehicle history and current profile, averaged over 20 independent trials and displayed in Table \ref{tab:results}.  Fig.\ref{fig:all} shows the panoply of visualizations produced by each trial, and Appendix F contains the visualizations for each trial, organized for side-to-side comparisons of different trials.
        
        \begin{table}[h]
            \centering
            \begin{tabular}{c|c|c|c|c|c}
                Method &  Variant & Navigation Error & Current Error  &Nav err (no final GPS) & Curr err (no final GPS)\\
                1 & Basic & 151 & 5.30e-2 & 318 & 7.85e-2 \\
                2 & Higher-order & {\bf 97} & {\bf 3.67e-2} & {\bf 216} & {\bf 5.46e-2} \\
                3 & Vehicle-Current covariance & 134 & 4.92e-2 & 308 & 7.40e-2 \\
                4 & Covariance and Higher-order & 120 &  6.56e-2 & 264 & 6.72e-2\\
            \end{tabular}
            \caption{Results, with the optimal method for each column boldfaced.  Navigation RMSE uses units of meters, current RMSE uses units of meters per second.  All methods perform competitively, with a different method taking top spot for each metric.}
            \label{tab:results}
        \end{table}
    
\begin{figure}
            \centering
            \includegraphics[width=35pc]{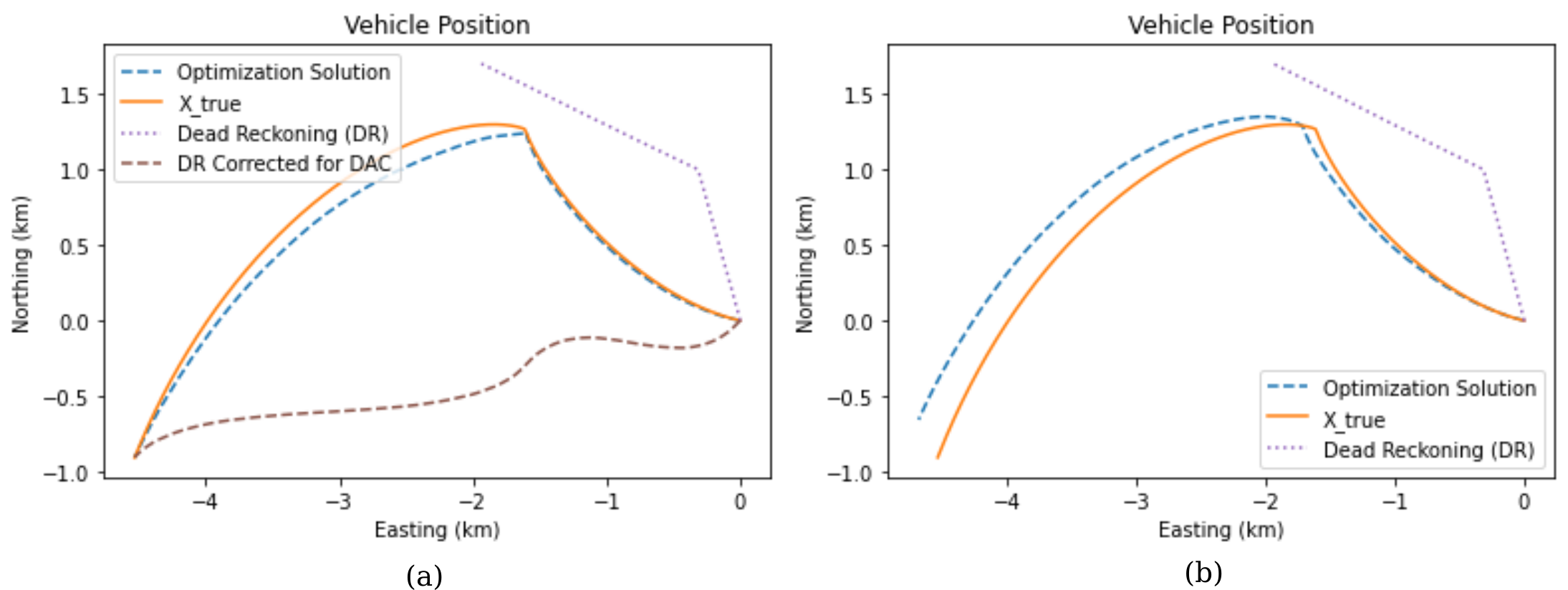}
            \caption{Navigation error (a) method 1, with start and end GPS points, and (b) method 4, including only starting GPS points.  All of this paper's methods perform well, outperforming DR corrected by DAC, and the best and track the vehicle to 55m/hr before receiving a resurfacing GPS fix.}
            \label{fig:nav-example}
        \end{figure}

\begin{figure}
            \centering
            \includegraphics[width=110mm]{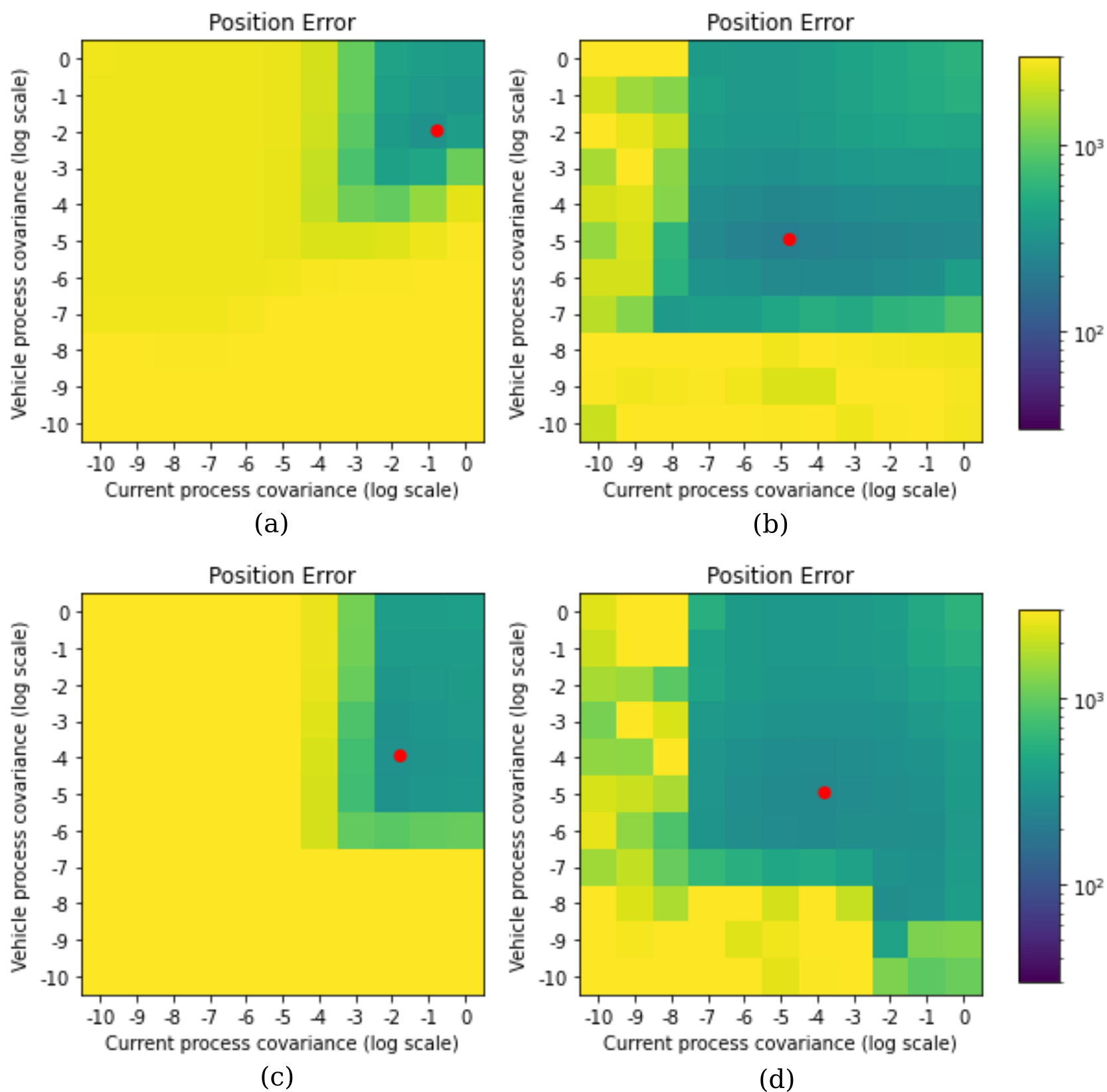}
            \caption{Navigation error across process variance parameter search (a) in method 1, (b) method 2, (c) method 3, and (d) method 4.  All cases shown include only starting GPS points, and the higher-order and correct-covariance smoothing methods have a much larger acceptable region of parameter space.}
            \label{fig:nav-search-example}
        \end{figure}

        All solution methods appear to perform adequately, with method 2 performing the best.  Considering first navigation error, most methods track the vehicle with at most a 150m difference between true and modeled position.  Even with just starting GPS positions, the second method keeps track of the vehicle with an accuracy of 80m after four hours in the test case.  For this particularly challenging track, the level of accuracy displayed in Fig.\ref{fig:nav-example} reflects an impressive success.  Looking in more detail at how the choice of vehicle and current process parameters affect navigation error in Fig.\ref{fig:nav-search-example}, we see that the higher order methods (2 and 4) and correct covariance methods (3 and 4) allow a much wider region of parameter space to produce satisfactory results.
        
\begin{figure}
            \centering
            \includegraphics[width=35pc]{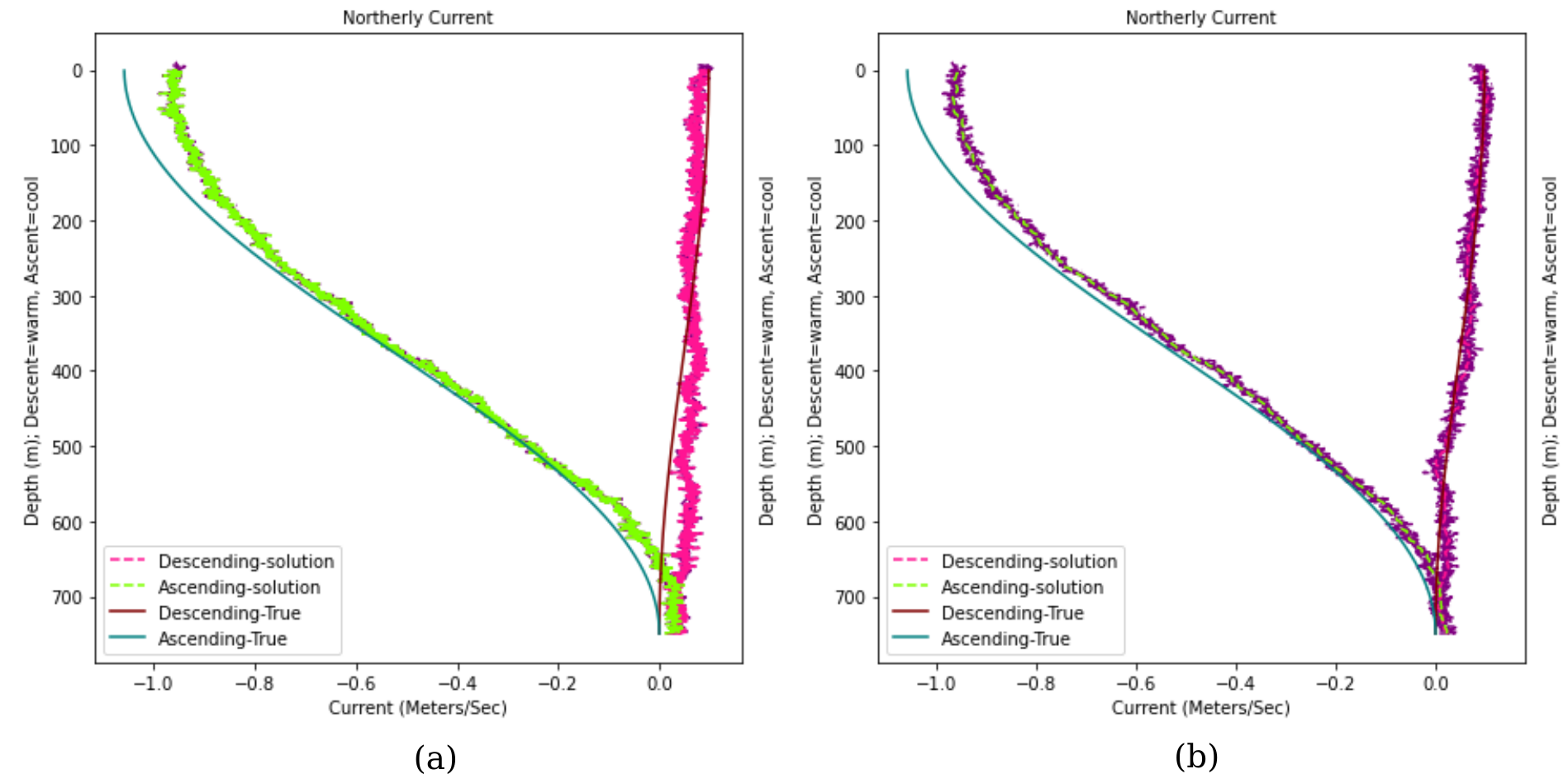}
            \caption{Current error (a) in method 3 and (b) method 4.  Both solutions shown include only starting GPS points.  Both the results on prepared metrics and the much smoother profile of the higher order methods recommend methods 2 and 4.}
            \label{fig:current-example}
        \end{figure}

        As for current error, most methods achieve around 6 cm/s error on average, with the strongest currents experiencing the highest error.  Given that our simulation includes several thousand depth points, this equates to a DAC error on the order of $10^{-5}$.  The higher order methods show a much less noisy solution; Figure \ref{fig:current-example} illustrates the comparison between methods 3 and 4 for cases including only starting GPS points.  The higher-order methods, methods 2 and 4, in the variant without a final GPS point both outperform the other methods and appear visually nicer.  Likewise, the much wider acceptable parameter region for higher order methods makes them better candidates for a production model.  We omit a separate plot for current error across the process covariance parameter search, as it looks similar to Fig.\ref{fig:nav-search-example}.  Among the higher-order methods, using the correct covariance provides a slightly smoother plot and larger acceptable parameter space, though averaged slightly worse performance.  
        
\begin{figure}
            \centering
            \includegraphics[width=16pc]{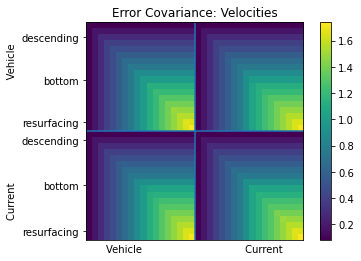}
            \caption{Method 1 velocity error covariance without a final GPS point appears appropriately similar to the covariance of Brownian motion.  Additionally, the repeated off-diagonal terms indicate a lack of identifiability.  All methods including only starting GPS points demonstrate these effects.}
            \label{fig:uq-brownian}
        \end{figure}

        The estimator variance plots tell several interesting stories.  Simulations including only beginning GPS points, such as in Fig.\ref{fig:uq-brownian}, show position and velocity uncertainty growing with duration since the initial GPS observations---as expected.  These square-corner covariance plots also appear to exactly replicate those of Brownian motion, providing more confidence that our model performs consistent with assumptions.

\begin{figure}
            \centering
            \includegraphics[width=35pc]{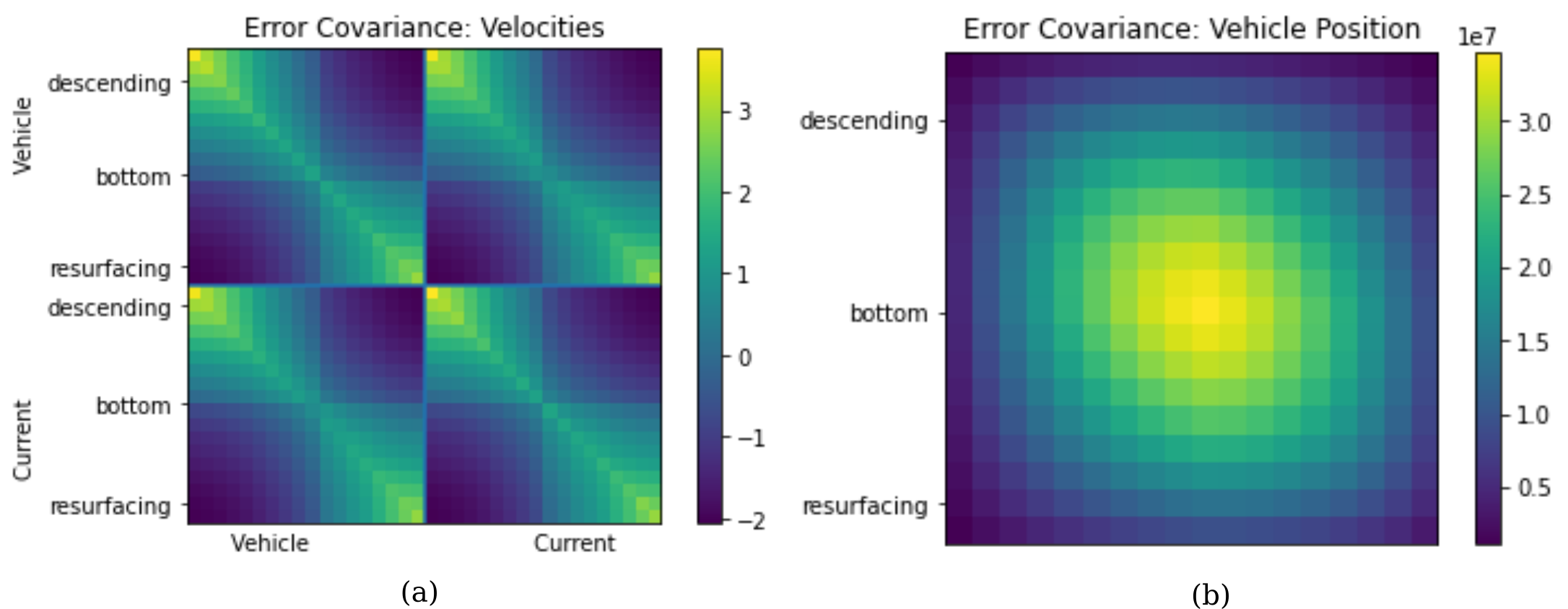}
            \caption{Method 1 (a) velocity and (b) position estimator error covariance including both starting and ending GPS measurements.  This example shows that position uncertainty reaches a minimum in vicinity of GPS observations, while velocity uncertainty reaches a maximum.}
            \label{fig:uq-pos-vs-vel}
        \end{figure}

        Simulations that included both start and finish GPS measurements, as in Fig.\ref{fig:uq-pos-vs-vel}, predictably show low position uncertainty at the endpoints and high uncertainty at nadir.  On the other hand, velocity uncertainty reaches maxima at the endpoints and a minimum at nadir.  Given that GPS measurements affect the integral of velocity, such terms constrain velocities most at the midpoint, where the velocity has the most adjacent points for smoothing.  This effect shows that additional measurements reduce uncertainty of the observed quantity, but the smoother determines how that information passes to different process terms.
    
        In all cases, the off-diagonal blocks show that vehicle velocity terms correlate highly with their simultaneous currents.  Such structure indicates a near rank-deficiency in our model, a lack of identifiability in whether an individual velocity reflects vehicle motion or currents.  It suggests that direct observation of a just a few true vehicle velocities or currents during the dive would dramatically improve estimates.  Indeed, \citet{Medagoda2016} finds that such a dramatic correction in navigation estimates occurs as soon as the vehicle gets a DVL lock on the seafloor.
        
    
        As a final note, while some of our second-order results may appear less smooth than those of other authors', much of that roughness merely reflects the decision not to bin observations.
        
\section{Conclusion and Further Research}
\label{conclusion}
    We demonstrated how the choice of smoother for current and vehicle motion reflects a prior belief in a probabilistic process and how that choice of belief affects the fit of the data.  We generated data to show that our basic model, essentially a method of \citet{Jonker2019} without binning, performed reasonably well.  Extending the model to higher order to enforce smoothness of velocities improved the model substantially and results in a much wider region to choose acceptable variance parameters.  Our other innovation, smoothing vehicle velocity and current velocity with proper covariance widened the region of acceptable smoothing parameters slightly and provided a more rigorous if not more successful model. However, it corrects the assumption that the current and vehicle motion are completely independent.  This correction would be more significant in larger currents and with sparser measurements.  All variants significantly outperformed dead reckoning.  In particular, models succeeded without the final GPS point as long as two nearby GPS observations before diving helped identify true surface drift. 
    
    Various ancillary corrections may help improve the evidence for our major contributions.  Most significantly, we would prefer to return to binning observations.  Matrix conditioning caused issues with higher order models and even with the proper-covariance models.  Additionally, measurement terms could benefit from reintroducing ADCP bias, compass error, and biofouling errors.  Preliminary tests of our model on real data shows that the sensors are usually biased, a more substantial separation of GPS points on the surface is required than just a few minutes, and putative ground truth measurements also suffer from inaccuracy.
    
    As a key avenue for improvement, one could replace the direct matrix inversion of our method with a generic convex solver (which produces equivalent results in the least-squares case).  The new solver, in turn, would allow for robust vehicle and current processes other than Brownian motion, a topic of interest in Kalman smoothing literature such as \cite{aravkin2017generalized}.  That paper notes robust process models could would better handle the idiosyncracies we found in real data, such as very strong ADCP outliers at apogee.  Other processes such as an Ornstein-Uhlenbeck/Gauss-Markov process could also help by introducing the assumption that current centers around zero and the vehicle velocity centers around its design drift velocity.  Appendix E has the likelihood expression for an OU process, which is amenable to the inversion solution method. 
    
    A convex solver would also permit more robust measurement forms than the linear least-squares case.  Thus, we could include nonlinear range measurements in our framework, as from Ultra-Short Baseline (USBL) or Long Baseline (LBL) systems or even the range-azimuth-elevation of the USBL method in \citet{Jakuba}.  For measurements close to $\vec c_x$ and $\vec c_y$, the local coordinates of the sensors, range measurement adds a log likelihood function of:
    \[
        -\ell(\vec z_{range}|x) = \frac{1}{2\sigma_{range}^2}\|\vec z_{range} - \left((\mat B^q_{range}\vec x_q - \vec c_x)^2 - (\mat B^q_{range}\vec x_q - \vec c_y)^2\right)^{1/2}\|^2,
    \]
    where the vector exponents apply elementwise.  These measurements, as nonlinear, give rise to a nonconvex likelihood expression.  Even so, the expression may still permit finding an acceptable local solution.
    
    Finally, automatic tuning of process covariance parameters could benefit usability of the model. While we obtained stable results in simulation using grid search for these parameters, a comprehensive approach for simultaneous state and current estimation with more efficient ways of finding underlying variance parameters would more automatically solve the engineering problem.
    
    The ability to calculate navigation and current state in real time without external submerged position support allows for more effective use of gliders.  For example, the vehicle could choose to navigate at a depth with beneficial currents to extend its range.  It could also dispose of the raw measurements and store the current profile in a compressed form, allowing longer dives and longer deployments.  In addition, the proper covariance smoothing could allow the solution to degrade less with infrequent measurements, saving battery life.  Such innovations would vastly improve the ability to understand global subsurface currents.

\appendix[A]
\appendixtitle{Kalman Matrices}
\label{appendix-kalman}
    The Kalman smoother arrives in our likelihood expression \eqref{eqn:basic-likelihood} in the form of a prior.  Let $W$ refer to a Brownian motion.  Referring to time index $j$ for a coordinate with velocity $x_{\dot q}^j= W_{t_j}$ and position $x_q^j = \int_0^{t_j}W_r dr$, expression \eqref{eqn:vehicle-process} is exactly the likelihood of Brownian motion and its integral.  
    That is: 
    \begin{align*}
        x^j_{\dot q} -x_{\dot q}^{j-1} &= W_{t_j}-W_{t_{j-1}}\\
        &\sim W_{\Delta t}\\
        x^j_q -x_q^{j-1} &= \int_{t_{j-1}}^{t_j} x_{\dot q}(t)dt\\
        &= \int_{t_{j-1}}^{t_j} x_{\dot q}^{j-1}+ x_{\dot q}(t) -x_{\dot q}^{j-1} dt\\
        &= \Delta t x_{\dot q}^{j-1} + \int_{t_{j-1}}^{t_j} W_t-W_{t_{j-1}}dt\\
        &\sim \Delta t x_{\dot q}^{j-1} + \int_0^{\Delta t}W_r dr.
    \end{align*}
    We have used the independence of successive increments and stationarity to change from $[t_{j-1}, t_j]$ to simpler variables $[0, \Delta t]$.  Thus,
    \begin{align}
    \label{eqn:kalman-var}
    \left[\begin{matrix}
            x^j_{\dot q} -x_{\dot q}^{j-1} \\
            x^j_q -x_q^{j-1} - \Delta t x_{\dot q}^{j-1}
        \end{matrix}\right]
        &\sim \left[\begin{matrix}
            W_{\Delta t} \\ \int_0^{\Delta t}W_r dr
        \end{matrix}\right]
    \end{align}
    Taken in blocks, matrix $\mat G_\upsilon \vec x_\upsilon$ is the left side of expression \eqref{eqn:kalman-var}.
    \begin{align}
	    \label{eqn:basic-kalman-update-appa}
    	\mat G_\upsilon = 
    	\left[\begin{matrix}
        	-1& 0 & 1 & 0 & &\\
        	-\Delta t_1 & -1 & 0 & 1 & &\\
        	& & -1 & 0 & 1& 0\\ 
        	& & -\Delta t_2 & -1 & 0 & 1\\
        	&&&&&& \ddots
    	\end{matrix}\right]~~~\text{and}~~~~
    	\vec x_\upsilon = \left[\begin{matrix}
	    x_{\dot q}^1\\
	    x_q^1\\
	    x_{\dot q}^2\\
	    x_q^2\\
	    \vdots
	    \end{matrix}\right].
	\end{align}
	Because the right hand side are integrals of Brownian motion, we know:
	\begin{align*}
	    \mathbb E[\mat G_\upsilon \vec x_\upsilon] = 0,
	\end{align*}
	With the mean determined, all that remains is to determine the covariance for each timestep: $\mat Q_\upsilon $ in the expression $\vec x^j_\upsilon|\vec x^{j-1}_\upsilon\sim \mathcal N((\mat G_\upsilon \vec x)^j, \mat Q_\upsilon^j)$.  Because increments are independent, the matrix $\mat Q_\upsilon$ has block diagonals $\mat Q_\upsilon^j$:
    \begin{align}
        \label{eqn:basic-kalman-covariance}
        \mathbb E\left[\begin{matrix}
            W_{\Delta t}^2 & W_{\Delta t}\cdot \int_0^{\Delta t}W_r dr \\
            W_{\Delta t}\cdot \int_0^{\Delta t}W_r dr & \left(\int_0^{\Delta t}W_r dr\right)^2
        \end{matrix}\right]
            =
    	\left[\begin{matrix}
    	\Delta t& \Delta t^2/2\\
    	\Delta t^2/2 & \Delta t^3/3
    	\end{matrix}\right].
    \end{align}
    We have omitted the calculation details involving Itô calculus.  In $\vec x_c$, where the state vector only includes a velocity, one can simply drop the rows/columns of $\mat G$ and $\mat Q$ corresponding to position.  Finally, since current varies by depth $s$, we replace $\Delta t$ with $\Delta s$.
    
\appendix[B]
\appendixtitle{Higher Order Smoothing}
\label{appendix-higher-order}
    In cases when we include an acceleration term, we can simply assign $x_{\ddot q}^j= W_{t_j}$ and  $x_{\dot q}^j = \int_0^{t_j}W_r dr$, giving the results from the previous appendix to the $x_{\ddot q}$ and $x_{\dot q}$ terms.  We need to integrate the Brownian motion one more time to get position:
    \begin{align*}
        x^j_q-x^{j-1}_q &=  \Delta t x_{\dot q}^{j-1}+\int_{t_{j-1}}^{t_j} x_{\dot q}(t) -x_{\dot q}^{j-1} dt \\
        &=\Delta t x_{\dot q}^{j-1}+\int_{t_{j-1}}^{t_j} \int_{t_{j-1}}^t x_{\ddot q}(r)dr dt\\
        &=\Delta t x_{\dot q}^{j-1}+\int_{t_{j-1}}^{t_j} \int_{t_{j-1}}^t x_{\ddot q}^{j-1}+x_{\ddot q}(r)-x_{\ddot q}^{j-1}dr dt\\
        &=\Delta t x_{\dot q}^{j-1}+\frac{\Delta t^2}{2}x_{\ddot q}^{j-1} +\int_{t_{j-1}}^{t_j} \int_{t_{j-1}}^t W_r-W_{t_{j-1}}dr dt\\
        &\sim\Delta t x_{\dot q}^{j-1}+\frac{\Delta t^2}{2}x_{\ddot q}^{j-1} +\int_0^{\Delta t} \int_0^t W_r dr dt\\
        &\sim\Delta t x_{\dot q}^{j-1}+\frac{\Delta t^2}{2}x_{\ddot q}^{j-1} +\int_0^{\Delta t} (\Delta t - t) W_r dt,\\
    \end{align*}
    where the last step is Cauchy's formula for iterated integrals.  Collecting terms as in the previous appendix, these relationships give rise to the formulae for $\mat G_\upsilon$ and $\mat Q_\upsilon$: 
	\begin{align}
    	\label{eqn:higher-order-kalman-appb}
    	\mat G_\upsilon = 
    	\left[\begin{matrix}
        	-1& 0 & 0 & 1 & 0 & 0& &\\
        	-\Delta t_1 & -1 & 0 & 0 & 1 & 0& & &\\
        	-\frac{\Delta {t_1}^2}{2}& -\Delta t_1 & -1 &0 & 0 & 1 & \\
        	& & & -1 & 0 & 0& 1 & 0 & 0 \\ 
        	& & & -\Delta t_2 & -1 & 0 & 0 & 1 & 0\\
        	& & & -\frac{\Delta {t_2}^2}{2}& -\Delta t_2 & -1 & 0 & 0 & 1\\
        	& & & & & & & & & \ddots
    	\end{matrix}\right],
    	\qquad 
    	 \mat Q_\upsilon=\left[\begin{matrix}
        	{\Delta t} & \frac{\Delta t^2}{2} & \frac{\Delta t^3}{6}\\
        	\frac{\Delta t^2}{2} & \frac{\Delta t^3}{3} & \frac{\Delta t^4}{8}\\
        	\frac{\Delta t^3}{6} & \frac{\Delta t^4}{8} & \frac{\Delta t^5}{20}\\
    	\end{matrix}\right].
	\end{align}

    The current process now has two total orders, and so has similar terms to vehicle smoothing in Appendix A, but again replacing $t$ with $s$.

\appendix[C]
\appendixtitle{Modeling Vehicle Relative Velocity}
\label{appendix-ttw}
    Return to equation \ref{change-of-variables}, where we sought to understand the implication of smoothing $x_r$, rather than $x_q$, and recognized $x_p$ terms.  If we focus on the second integral,
    	\[
    		\int_{s(t_1)}^{s(t_2)}x_{\dot p}(s, t(s))\frac{1}{\dot s} ds,
    	\]
        we ought to make several assumptions explicit:
        \begin{enumerate}
            \item Current does not change in the time or horizontal strata between observations
            \item Depth rate remains constant between two nearby observations.
        \end{enumerate}
        The first assumption simplifies our $x_{\dot p}$ term, and the second allows the $1/\dot s$ becomes a coefficient:
        \begin{align*}
            =\frac{1}{\dot s}\int_{s(t_1)}^{s(t_2)}x_{\dot p}(s) ds = \frac{\Delta s \cdot x_{\dot p}(t_1)+\int_{s(t_1)}^{s(t_2)}dW_r dr}{\dot s},
        \end{align*}
        Thus, our current process term, equation \ref{current_process}, now uses the process update matrix:
    	\begin{align}
    	    \label{eqn:ttw-kalman-update}
        	\mat G_c = 
        	\left[\begin{matrix}
            	-1& 0 & 1 & 0 & &\\
            	-\Delta s_1/\dot s_1 & -1 & 0 & 1 & 0&\\
            	0 & 0 & -1 & 0 & 1& \ddots\\ 
            	0 & 0 & -\Delta s_2/\dot s_2 & -1 & \ddots & \ddots
        	\end{matrix}\right],~~~~~~~
        	\vec x_c = \left[\begin{matrix}
    	    x_{\dot p}^1\\
    	    x_p^1\\
    	    \vdots\\
    	    x_{\dot p}^n\\
    	    x_p^n
    	    \end{matrix}\right].\\
    	\end{align}
        Likewise, $\mat Q_c$ now has blocks:
        \begin{align}
            \label{eqn:ttw-kalman-covariance}
    	    \left[\begin{matrix}
    	    \Delta s &  \Delta s^2/(2\dot s)\\
    	    \Delta s^2/(2\dot s) & \Delta s^3/(3\dot s^2)
    	    \end{matrix}\right].~~~~~~~~~~~
        \end{align}

\appendix[D]
\appendixtitle{Proper Covariance for Current and Vehicle True Velocity}
\label{appendix-covariance}
    We seek a specification for $\mathrm{Pr}(\vec x_{\dot q},~\vec x_q|~\vec x_c)$ in equation \eqref{eqn-conditional}.  We first form the correct covariance matrix for $\mathrm{Pr}(\vec x_{\dot q},~\vec x_q,~\vec x_c)$ and then derive the conditional covariance matrix and expected value.  Since $\mat G_\upsilon$ uses the expected value to recenter the process, we need only concern ourselves with a single timestep $\Delta t$.
    
    Although we do not utilize an $\vec x_p$ or $\vec x_r$ term, we still return to equation \eqref{relative}.  Writing out all the terms as random variables driven by independent Brownian motions $W^1_s$ driving current and $W^2_t$ driving the vehicle's speed through the water, we calculate:
    \begin{align*}
        x_c^j-x_c^{j-1} &= \int_{s(t_{j-1})}^{s(t_j)} \sigma_c^2 dW^1_s \sim \int_{0}^{\Delta s} \sigma_c^2 dW^1_s\\
        x_{\dot q}^j - x_{\dot q}^{j-1} &= \int_{t_{j-1}}^{t_j} \sigma_\upsilon^2 dW^2_t+\int_{s(t_{j-1})}^{s(t_j)} \sigma_c^2 dW^1_s \sim \int_{0}^{\Delta t} \sigma^2_c dW^2_t + \int_{0}^{\Delta s} \sigma^2_\upsilon dW^1_s\\
        x_{q}^j - x_{q}^{j-1} &= \int_{t_{j-1}}^{t_j} dt \cdot x_{\dot q}(t)\\
        &= \Delta t \cdot x_{\dot q}^{j-1}+ \int_{t_{j-1}}^{t_j} dt \cdot \left(x_{\dot q}(t)- x_{\dot q}^{j-1}\right)\\
        \int_{t_{j-1}}^{t_j} dt \cdot \left(x_{\dot q}(t)- x_{\dot q}^{j-1}\right)&\sim \int_{t_{j-1}}^{t_j} dt \int_{t_{j-1}}^{t} \sigma^2_\upsilon dW^2_r+\int_{s(t_{j-1})}^{s(t)} \sigma^2_c dW^1_r\\
        &= \int_0^{\Delta t}dt\int_{0}^{t} \sigma^2_\upsilon dW^2_r + \int_0^{\Delta t}dt\int_{0}^{s(t)} \sigma^2_c dW^1_r\\
        &= \int_0^{\Delta t}dt\int_{0}^{t}\sigma^2_\upsilon dW^2_r + \frac{1}{\dot s}\int_0^{\Delta s}ds\int_{0}^{s} \sigma^2_c dW^1_r~~~\text{(since $dt/ds=\dot s^{-1}$)}.
    \end{align*}
    For further simplification, we now let $\Delta x_c = x_c^j-x_c^{j-1}$, $\Delta x_{\dot q} = x_{\dot q}^j - x_{\dot q}^{j-1}$, and $\Delta x_q = x_{q}^j - x_{q}^{j-1} -\Delta t \cdot x_{\dot q}^{j-1}$.  These mean-zero Gaussian random variables have the covariance matrix:
    \begin{align*}
        \mat Q = &\left[\begin{matrix}
        \sigma^2_c\Delta s & \sigma^2_c\Delta s & 
            \sigma^2_c\frac{\Delta s^2}{2\dot s} \\
        \sigma^2_c\Delta s & \sigma^2_\upsilon\Delta t + \sigma^2_c\Delta s &
            \sigma^2_\upsilon\frac{\Delta t^2}{2}+ \sigma^2_c\frac{\Delta s^2}{2\dot s}\\
        \sigma^2_c\frac{\Delta s^2}{2\dot s} &
            \sigma^2_\upsilon\frac{\Delta t^2}{2}+ \sigma^2_c\frac{\Delta s^2}{2\dot s} &
            \sigma^2_\upsilon\frac{\Delta t^3}{3}+\sigma^2_c\frac{\Delta s^3}{3\dot s^2} \\
        \end{matrix}\right].
        ~~\gray{\left. \begin{matrix}
        \Delta x_c\\\Delta x_{\dot q}\\ \Delta x_q
        \end{matrix}\right.} \\
        & \gray{\begin{matrix}~~~~~\Delta x_c~~~~~~~~~~~&\Delta x_{\dot q}~~~~~~~~~~~~~~~~~& \Delta x_q\end{matrix}}
    \end{align*}
    Calculating the conditional distribution gives:
    \begin{align}
    \label{eqn:conditional-rv}
        \Delta x_c& \sim \mathcal N(0, \Delta s)\\
        \left.\left[\begin{matrix}
            \Delta x_{\dot q}\\
            \Delta x_q
        \end{matrix}\right]\right|\Delta x_c
            &\sim \mathcal N \left(
                \left[\begin{matrix}
                    \Delta x_c \\
                    \frac{\Delta s}{2\dot s}\Delta x_c
                \end{matrix}\right],
                \sigma^2_\upsilon\left[\begin{matrix}
                    \Delta t  & \frac{\Delta t^2}{2}\\
                    \frac{\Delta t^2}{2} & \frac{\Delta t^3}{3}+\frac{\sigma^2_c}{\sigma^2_\upsilon}\frac{\Delta s^3}{12\dot s^2}
                \end{matrix}\right]
            \right).
    \end{align}
    When there are $k-1$ intermediate depths in the state vector between $s(t_{j-1})$ and $s(t_j)$, the terms with $\Delta s$ need to be summed across all $i\in \{1,\dots,k\}$, e.g. $\frac{\Delta s}{2\dot s}\Delta x_c \rightarrow \sum_{i=1}^k \frac{\Delta s_i}{2\dot s}\Delta x_{c_i}$ and $\frac{\Delta s^3}{12\dot s^2}\rightarrow \sum_{i=1}^k \frac{\Delta s_i^3}{12\dot s^2}$.  Depth rate $\dot s$ is assumed constant between timepoints.
    
    When smoothing the vehicle process to higher order and modeling acceleration conditionally, we also choose to smooth the current process to the same order.  Thus, to smooth current to higher order, we re-notate $\vec x_c$ as $\vec x_{\dot p}$, the current velocity, and $\vec x_{\ddot p}$, the change in current with depth.  Since $\vec x_{\ddot p}$ is a change in velocity with respect to depth, we must multiply by $\dot s$ in order to get the change in velocity with respect to time.  That is,
    \begin{align*}
        x_{\ddot q}^j - x_{\ddot q}^{j-1} &= \int_{t_{j-1}}^{t_j} \sigma_\upsilon^2 dW^2_t+\int_{s(t_{j-1})}^{s(t_j)} \dot s \sigma_c^2 dW^1_s.
    \end{align*}
    The $\Delta$ variables (e.g., $\Delta x_q$, $\Delta x_{\dot q}$) then also must adjust to maintain a mean-zero process in the joint distribution. The conditional distribution becomes:
    \begin{align}
    \label{eqn:conditional-rv-high-order}
        \left.\left[\begin{matrix}
            \Delta x_{\ddot q}\\
            \Delta x_{\dot q}\\
            \Delta x_q
        \end{matrix}\right]\right|
        \left[\begin{matrix}\Delta x_{\ddot p}\\\Delta x_{\dot p}\end{matrix}\right]
        \sim \mathcal N \left(
            \left[\begin{matrix}
                \dot s \Delta x_{\ddot p}\\
                \Delta x_{\dot p}\\
                \frac{\Delta s}{2\dot s}\Delta x_{\dot p} - \frac{\Delta s^2 \Delta x_{\ddot p}}{12\dot s}\\
            \end{matrix}\right],
            \sigma^2_\upsilon\left[\begin{matrix}
                \Delta t & \frac{\Delta t^2}{2} & \frac{\Delta t^3}{6}\\
                \frac{\Delta t^2}{2} & \frac{\Delta t^3}{3} & \frac{\Delta t^4}{8}\\
                \frac{\Delta t^3}{6} & \frac{\Delta t^4}{8} & 
                    \frac{\Delta t^5}{20}+\frac{\sigma^2_c}{\sigma^2_\upsilon}\frac{\Delta s^5}{720\dot s^2}\\
            \end{matrix}\right]
        \right).
    \end{align}
    In both expressions \ref{eqn:conditional-rv} and \ref{eqn:conditional-rv-high-order}, we create $\widetilde{\mat G}_\upsilon$ to subtract the conditional mean in order to recover a mean-zero random variable.

\appendix[E]
The Orstein-Uhlenbeck process is Gauss Markov with a constant drift rate $k$ and shifted by $\theta$:
\begin{equation}
    d\dot x_t = k(\theta - \dot x_t)dt + \alpha dW_t.
\end{equation}
Thus, the process always seeks to revert towards the value $\theta$, which could represent a prior average current or designed vehicle drift velocity.  We can calculate that the vehicle state at any time conditioned on its previous state is:
\begin{equation}
    \left.\left[\begin{matrix}
        x_{\dot q}^{j+1} \\ x_q^{j+1}
        \end{matrix}\right]\right|
        \left[\begin{matrix}
        x_{\dot q}^j \\ x_q^j
    \end{matrix}\right]
    \sim
    \mathcal N\left(
        \left[\begin{matrix}
            e^{-k\Delta t}x_{\dot q}^j+(1-e^{-k\Delta t})\theta\\
            x_q^j -\frac{1}{k}\left(1-e^{-k\Delta t}\right)x_{\dot q}^j +\left(\frac{1}{k}+\Delta t-\frac{e^{-k\Delta t}}{k}\right)\theta 
        \end{matrix}\right],
        \mat Q^j
    \right),
\end{equation}
with the covariance matrix,
\begin{equation}
    \mat Q^j = 
        \left[\begin{matrix}
            \frac{\alpha^2}{2k}\left(1-e^{-2k\Delta t}\right) &
                \frac{\alpha^2}{k^2}e^{-k\Delta t}\left(\cosh{k\Delta t}-1\right)\\
            \frac{\alpha^2}{k^2}e^{-k\Delta t}\left(\cosh{k\Delta t}-1\right) &
                \frac{\alpha^2}{k^2}\Delta t 
                    + \frac{\alpha^2}{k^3}\left(
                        2e^{-k\Delta t}-\frac{1}{2}e^{-2k\Delta t}-3/2
                    \right)
        \end{matrix}\right].
\end{equation}

This result allows us to use $\|\mat G\vec x-\vec b\|^2_{\mat Q^{-1}}$ as our prior negative log-likelihood with the same block diagonal structure as the original Kalman filter for a Brownian velocity.  Similar to the original Kalman filter, the OU Kalman filter can be extended to higher order or made conditional upon another process as in the previous appendices.

\appendix[F]
\appendixtitle{Complete Plots for Comparison}
\label{appendix-plots}
    \begin{figure}[h]
        \centering
        \includegraphics[width=120mm]{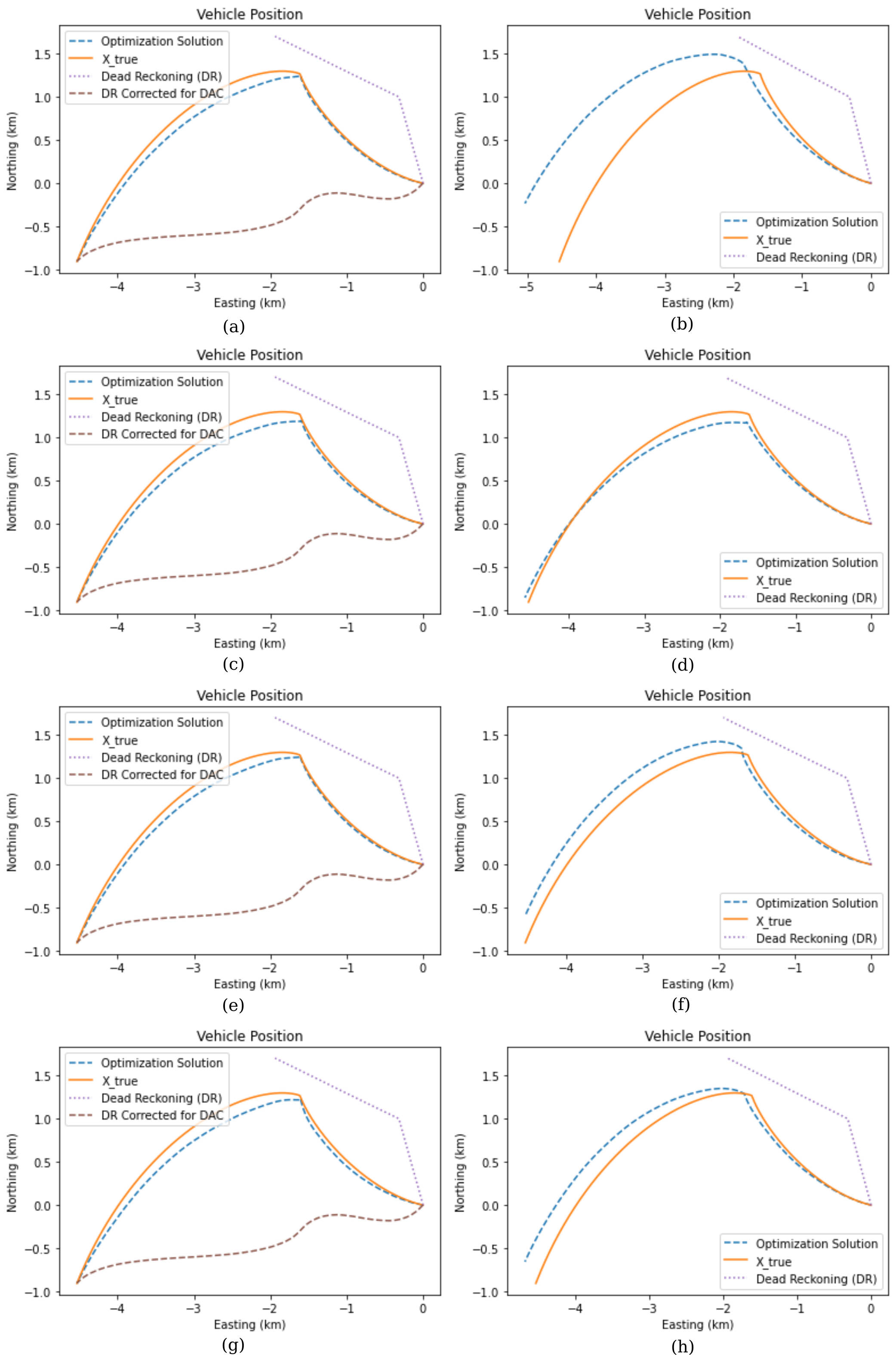}
        \caption{Representative navigation results from (a, b) method 1, (c, d) method 2, (e, f) method 3, and (g, h) method 4.  Trials on the left include start and end GPS points, trials on the right include two starting GPS points only}
        \label{fig:nav}
    \end{figure}
    \clearpage
    \begin{figure}[h]
        \centering
        \includegraphics[width=100mm]{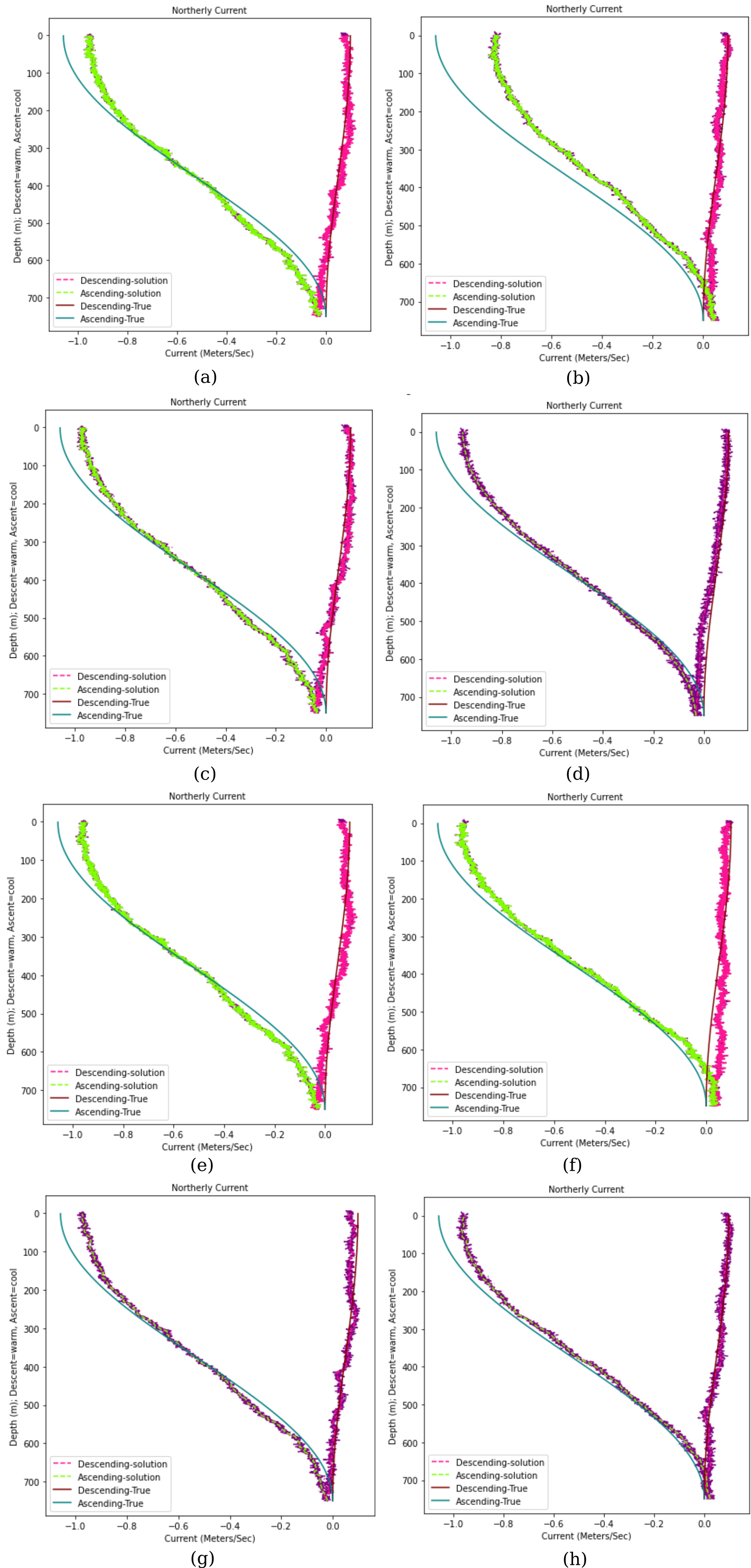}
        \caption{Representative current results from (a, b) method 1, (c, d) method 2, (e, f) method 3, and (g, h) method 4.  Trials on the left include start and end GPS points, trials on the right include two starting GPS points only}
        \label{fig:current}
    \end{figure}
    \clearpage
    \begin{figure}[h]
        \centering
        \includegraphics[width=110mm]{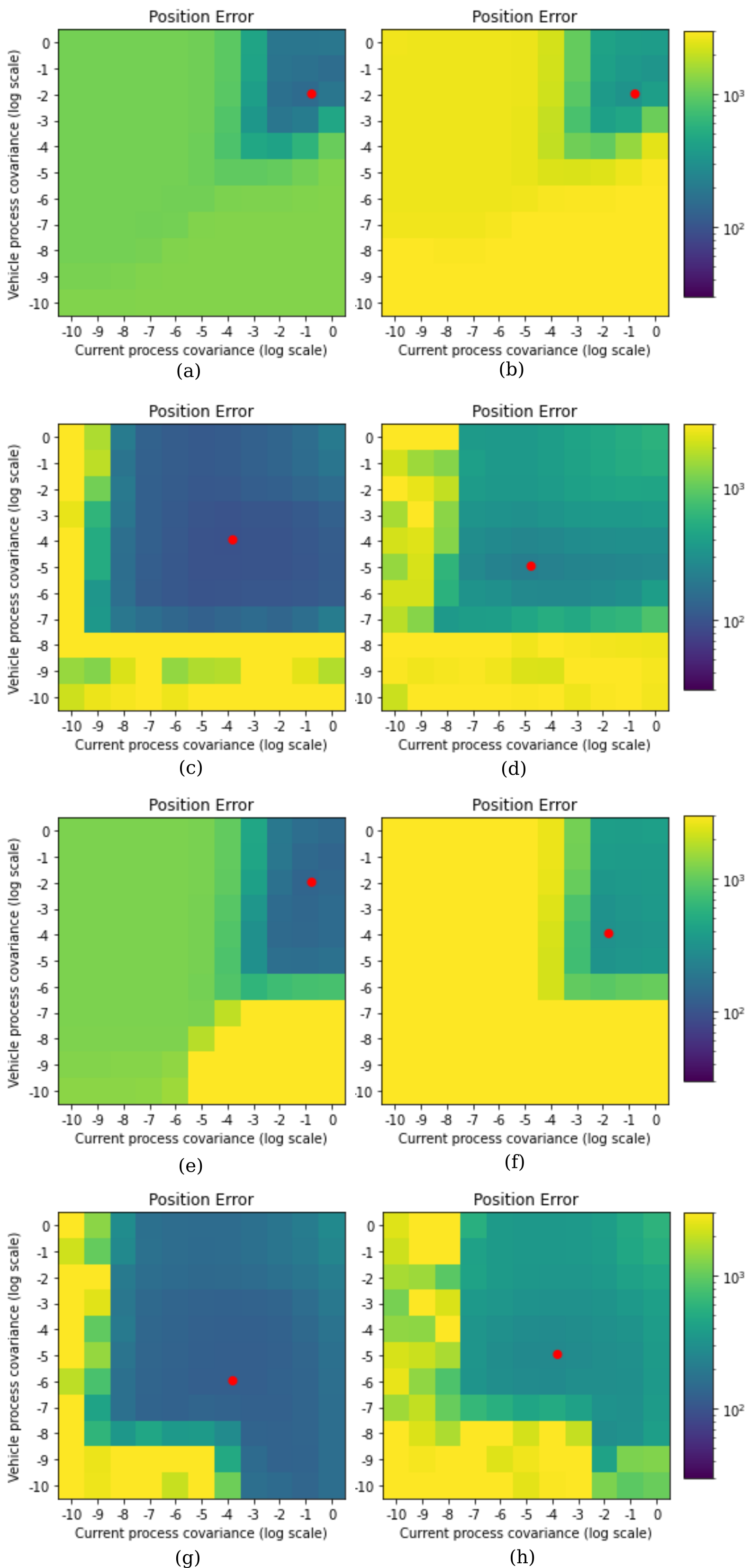}
        \caption{Navigation parameter search results from (a, b) method 1, (c, d) method 2, (e, f) method 3, and (g, h) method 4.  Trials on the left include start and end GPS points, trials on the right include two starting GPS points only.  All colors normalized to same values.}
        \label{fig:nav-search}
    \end{figure}
    \clearpage
    \begin{figure}[h]
        \centering
        \includegraphics[width=105mm]{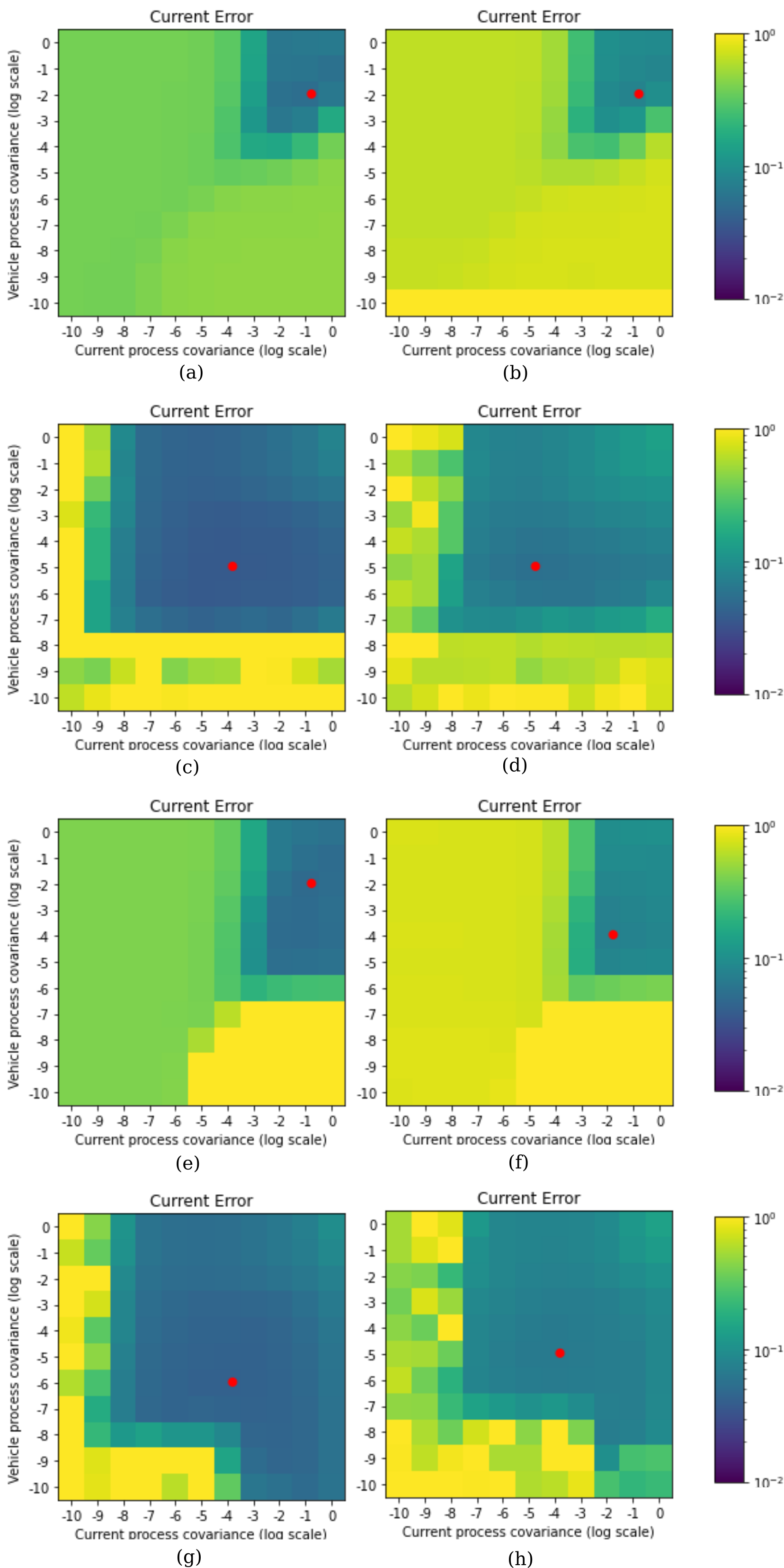}
        \caption{Current parameter search results from (a, b) method 1, (c, d) method 2, (e, f) method 3, and (g, h) method 4.  Trials on the left include start and end GPS points, trials on the right include two starting GPS points only.  All colors normalized to same values.}
        \label{fig:current-search}
    \end{figure}
    \clearpage
    \begin{figure}[h]
        \centering
        \includegraphics[width=105mm]{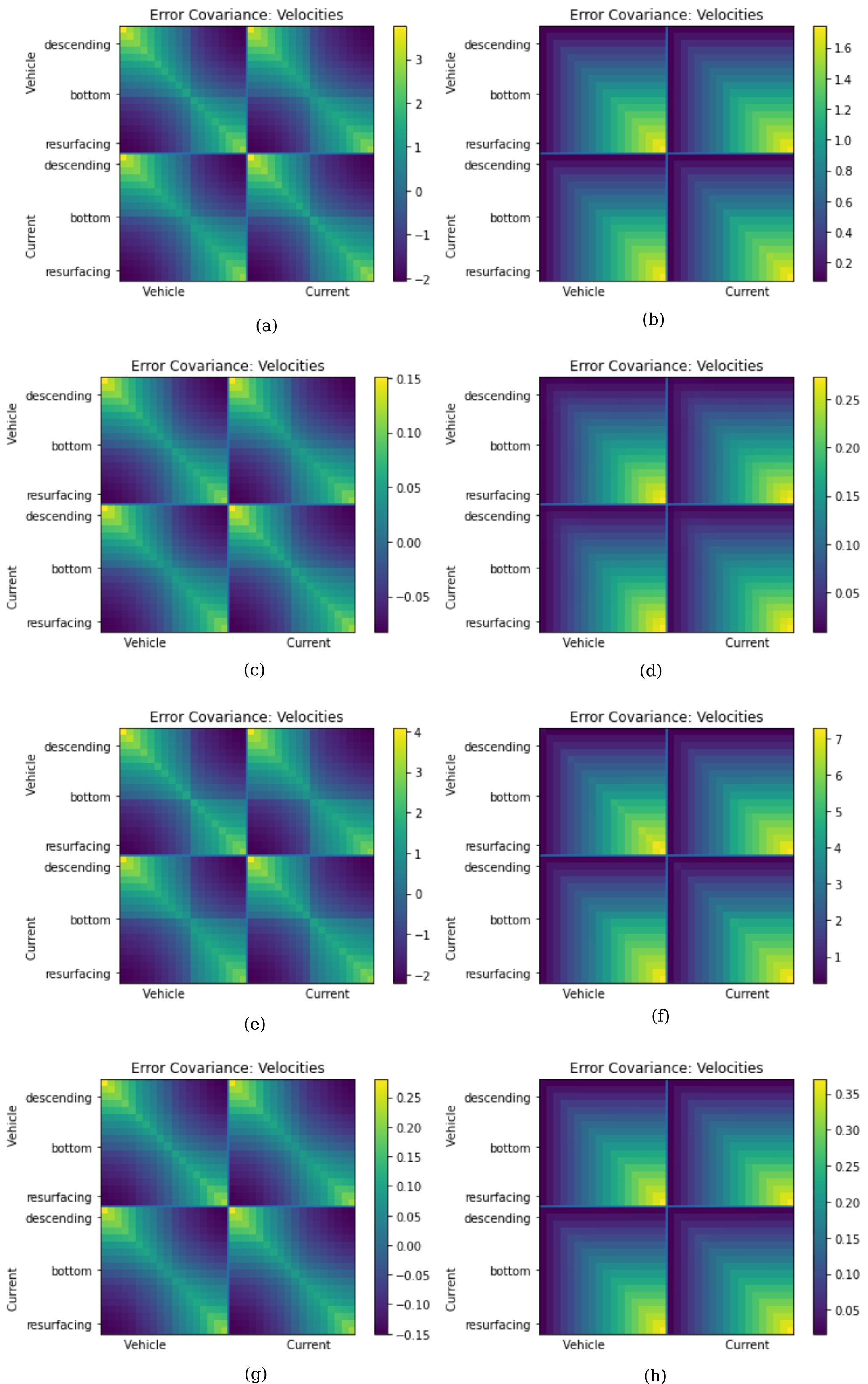}
        \caption{Velocity estimator covariance (subset of all states) from (a, b) method 1, (c, d) method 2, (e, f) method 3, and (g, h) method 4.  Trials on the left include start and end GPS points, trials on the right include two starting GPS points only.  Colors are {\bf not} normalized to same values.}
        \label{fig:uq}
    \end{figure}
    \clearpage
    \begin{figure}[h]
        \centering
        \includegraphics[width=105mm]{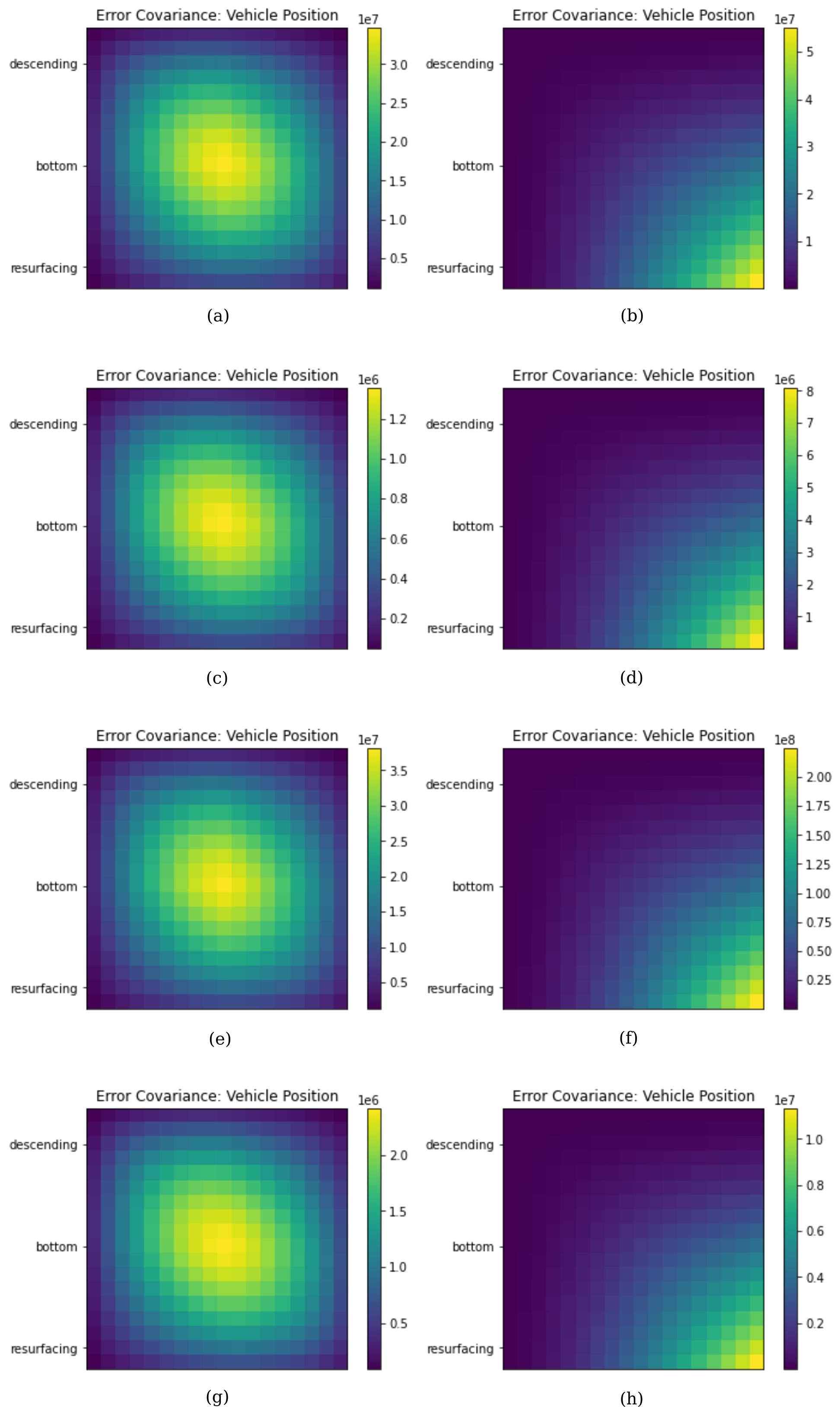}
        \caption{Navigation position estimator covariance (subset of all states) from (a, b) method 1, (c, d) method 2, (e, f) method 3, and (g, h) method 4.  Trials on the left include start and end GPS points, trials on the right include two starting GPS points only.  Colors are {\bf not} normalized to same values.}
        \label{fig:uq-posit}
    \end{figure}

\acknowledgments
This research was supported by the NOAA Office of Exploration and Research award NA20OAR0110429 and by the Veterans Administration under the GI Bill.

%
%
\datastatement
One can generate the figures and data for this paper by downloading our Python package at https://github.com/UW-AMO/seaglider-navigation and running "publish\_figures.py".  The commit used to generate the exact figures in the paper is e63028d.

%






%



\bibliography{adcp}

\end{document}